\documentclass[11pt,a4paper]{amsart}

\usepackage{amsfonts, amsmath, amssymb, amsthm}
\usepackage{amsmath, amssymb}

 \theoremstyle{plain}
 \newtheorem{thm}{Theorem}[section]
 \newtheorem*{thm*}{Theorem}
 \newtheorem{lem}[thm]{Lemma}
 \newtheorem{prop}[thm]{Proposition}
 \newtheorem{cor}[thm]{Corollary}
 \theoremstyle{definition}
 
 \newtheorem{exmp}[thm]{Example}

 \theoremstyle{remark}
 \newtheorem{rem}[thm]{Remark}
 
 \newtheorem{prob}[thm]{Problem}
 
\usepackage{txfonts}
\usepackage{mathbbol}
\usepackage{ae}

\usepackage{color}
\usepackage{graphicx}
\usepackage{url}
\usepackage{overpic}
\usepackage{subfigure}
\usepackage{booktabs}
\usepackage{psfrag}
\usepackage{paralist}
\usepackage[a4paper]{anysize}
\usepackage[all,ps]{xy}

\newcommand\cR{{\mathcal R}}

\newcommand\NN{{\mathbb N}}
\newcommand\RR{{\mathbb R}}
\newcommand\QQ{{\mathbb Q}}
\newcommand\ZZ{{\mathbb Z}}
\newcommand\CC{{\mathbb C}}
\newcommand\PP{{\mathbb P}}
\newcommand\Sph{{\mathbb S}}
\newcommand\tropPG{{\mathbb T\mathbb P}}
\newcommand\abs[1]{\left|#1\right|}

\newcommand\SetOf[2]{\left\{#1\,\vphantom{#2}\right.\left|\vphantom{#1}\,#2\right\}}
\newcommand\smallSetOf[2]{\{#1\,|\,#2\}}
\newcommand\conv{\operatorname{conv}}

\newcommand\stc{\operatorname{stc}}

\newcommand\vol[1]{\operatorname{vol}({#1})}
\newcommand\simpp{\times_{\operatorname{stc}}}
\newcommand\Initial{\operatorname{in}}
\newcommand\restr{\kern-1ex\mid}

\renewcommand\mod{\operatorname{mod}}
\renewcommand\phi{\varphi}

\renewcommand\phi{\varphi}

\graphicspath{{images/}}

\begin{document}

\title{Products of Foldable Triangulations}
\author[Joswig and Witte]{Michael Joswig \and Nikolaus Witte}
\address{Michael Joswig, Fachbereich Mathematik, AG~7, TU Darmstadt, 64289 Darmstadt, Germany}
\email{joswig@mathematik.tu-darmstadt.de}
\address{Nikolaus Witte, Fachbereich Mathematik, AG~7, TU Darmstadt, 64289 Darmstadt, Germany}
\email{witte@math.tu-berlin.de}
\thanks{Both authors are supported by Deutsche Forschungsgemeinschaft, DFG Research Group
  ``Polyhedral Surfaces.''}
\date{\today}

\begin{abstract}
  Regular triangulations of products of lattice polytopes are constructed with the additional
  property that the dual graphs of the triangulations are bipartite.  The (weighted) size difference
  of this bipartition is a lower bound for the number of real roots of certain sparse polynomial
  systems by recent results of Soprunova and Sottile [Adv.\ Math.\ 204(1):116--151, 2006].
  Special attention is paid to the cube case.
\end{abstract}

%

\maketitle

\section{Introduction}

A triangulation $K$ of an $m$-polytope $P$ is \emph{foldable} if $K$ admits a non-degenerate
simplicial map to an $m$-simplex.  This is equivalent to the property that its $1$-skeleton is
colorable in the graph-theoretic sense with the minimally possible number of $m+1$ colors. Further,
a triangulation is \emph{regular} if it can be lifted to $m+1$ dimensions as a lower convex hull.
The barycentric subdivision of any regular triangulation is an example of a triangulation which is
both regular and foldable.  A lattice triangulation of $P$ is \emph{dense} if its vertices are all
the lattice points inside~$P$, and, for the sake of brevity, we refer to a regular, dense, and
foldable triangulation as an \emph{{rdf}-triangulation}.  It is known that a triangulation of a
polytope (or, more generally, any simply connected manifold) is foldable if and only if its dual
graph is bipartite; see~\cite{MR1900311}.  From {rdf}-triangulations of lattice polytopes Soprunova
and Sottile~\cite{math.AG/0409504} construct sparse polynomial systems with non-trivial lower bounds
for the number of real roots.

For generic coefficients the exact number of complex solutions of a sparse system of polynomials is
known from Kushnirenko's Theorem~\cite{kushnirenko:newton_polyhedron}.  To estimate the number of
real solutions, however, is considerably more delicate.  The lower bound in the approach of
Soprunova and Sottile is the degree of a map on the oriented double cover of the real part $Y_P$ of
the toric variety associated with the lattice polytope $P$, where $P$ comes in as the common Newton
polytope of the polynomials in the system.  In combinatorial terms this map degree translates into
the size difference of the two color classes of facets of a
{rdf}-triangulation $K$ of~$P$.  More precisely, only those facets of~$K$ count in the size difference,
called the \emph{signature}, which have odd normalized volume.  We sketch this approach in
Section~\ref{sec:real_roots}.

This paper is mainly focused on the combinatorial aspects, but we apply our results to sparse
polynomial systems at the end.  We form {rdf}-triangulations of products of lattice polytopes
from {rdf}-triangulations of the factors.  As an application we construct triangulations of the
$d$-cube $C_d=[0,1]^d$, which is the product of~$d$ line segments.  Here we find
{rdf}-triangulations of~$C_d$ with a super exponentially large signature.  Optimizing triangulations of cubes for
combinatorial parameters is often difficult, and basic questions are still open: Most prominently,
for the minimal number of facets in a $d$-cube triangulation for $d>7$ only partial asymptotic
results are known; see Anderson and Hughes~\cite{MR1411113}, Smith~\cite{MR1737333}, Orden and
Santos~\cite{MR2013970}, Bliss and Su~\cite{math.CO/0310142}, and Zong~\cite{MR2133310}.  The
question whether the constructed triangulations of the $d$-cube have maximal signature is not
addressed in this paper.

The paper is organized as follows.  We start out with studying products of simplices because these
naturally form the building blocks in our product triangulations.  The key player here is the
staircase triangulation studied by Billera, Cushman, and Sanders~\cite{MR976522}, Gel$'$fand,
Kapranov, and Zelevinsky~\cite{MR1264417}, and others.  Then we focus on products of arbitrary
simplicial complexes.  These simplicial products, which depend on linear orderings of the vertices
of the factors, already occur in the work of Eilenberg and Steenrod~\cite[Section~II.8]{MR0050886}
and Santos~\cite{MR1758756}.  We prove that the product of two foldable simplicial complexes again
has a foldable triangulation.  Here it is important that there are still some choices left, a fact
which plays a role in the construction of the cube triangulations.  Then we can prove the following
Combinatorial Product Theorem, which is Theorem~\ref{thm:comb_product_thm} in this paper.

\begin{thm*}
  Let $P^\lambda$ and $Q^\mu$ be {rdf}-triangulations of an
  $m$-dimensional lattice polytope $P\subset\RR^m$ and an $n$-dimensional lattice polytope
  $Q\subset\RR^n$, respectively.  For specific vertex orderings of the factors (to be explained
  later) the simplicial product $P^\lambda\simpp Q^\mu$ is an {rdf}-triangulation of the
  polytope $P\times Q$ with signature
  \[\sigma(P^\lambda\simpp Q^\mu)=\sigma_{m,n}\;\sigma(P^\lambda)\;\sigma(Q^\mu)\;,\]
  where $\sigma_{m,n}$ is the signature of the staircase triangulation of the product of simplices
  $\Delta_m\times\Delta_n$.
\end{thm*}

For the algebraic applications it is essential that Theorem~\ref{thm:comb_product_thm} can
further be improved.  In Theorem~\ref{thm:alg_product_theorem} we show that (with a mild additional
assumption) the simplicial product $P^\lambda\simpp Q^\mu$ meets the geometric requirements of Soprunova and
Sottile, provided that both factors do.

As an application of our Product Theorems the paper continues with an explicit construction of
{rdf}-triangulations of the $d$-cube with signature in $\Omega(\lfloor d/2\rfloor!)$.
This lower bound partially relies on computational results obtained with
\texttt{TOPCOM}~\cite{topcom}, \texttt{polymake}~\cite{polymake,DMV:polymake,math.CO/0507273},
\texttt{MAGMA}~\cite{magma}, and \texttt{QEPCAD}~\cite{qepcad}.

\section{Products of Simplices}
\label{sec:prod_of_simp}

Let $\Delta_m=\conv(0,e_1,\dots,e_m)$ be the standard $m$-simplex, where $e_i$ denotes the $i$-th
unit vector of $\RR^m$.  Its normalized volume $\nu(\Delta_m)$ equals $\vol{\Delta_m}\,m!=1$.

The product $\Delta_m\times\Delta_n$ is an $(m+n)$-dimensional convex polytope with $(m+1)(n+1)$
vertices and $m+n+2$ facets.  As one key feature $\Delta_m\times\Delta_n$ has the property that it
is \emph{totally unimodular}, that is, each facet of any triangulation which uses no additional
vertices has normalized volume~$1$.
As a consequence the size of an arbitrary such triangulation of $\Delta_m\times\Delta_n$ is
\[\nu(\Delta_m\times\Delta_n)=\vol{\Delta_m}\vol{\Delta_n}\,(m+n)!=\binom{m+n}{m}\;.\]

\begin{figure}[bhtp]\centering
  \raisebox{.044\textwidth}{\includegraphics[height=.36\textwidth]{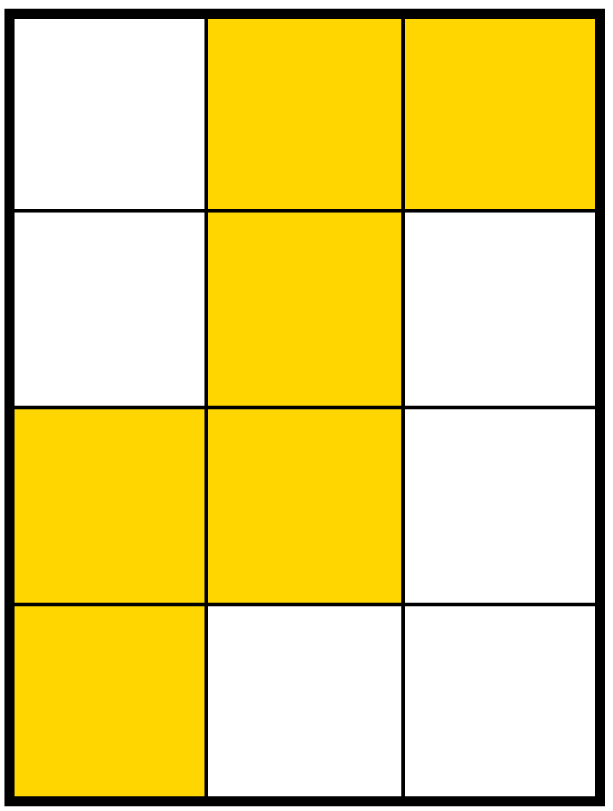}}
  \hskip.1\textwidth
  \begin{overpic}[width=.5\textwidth]{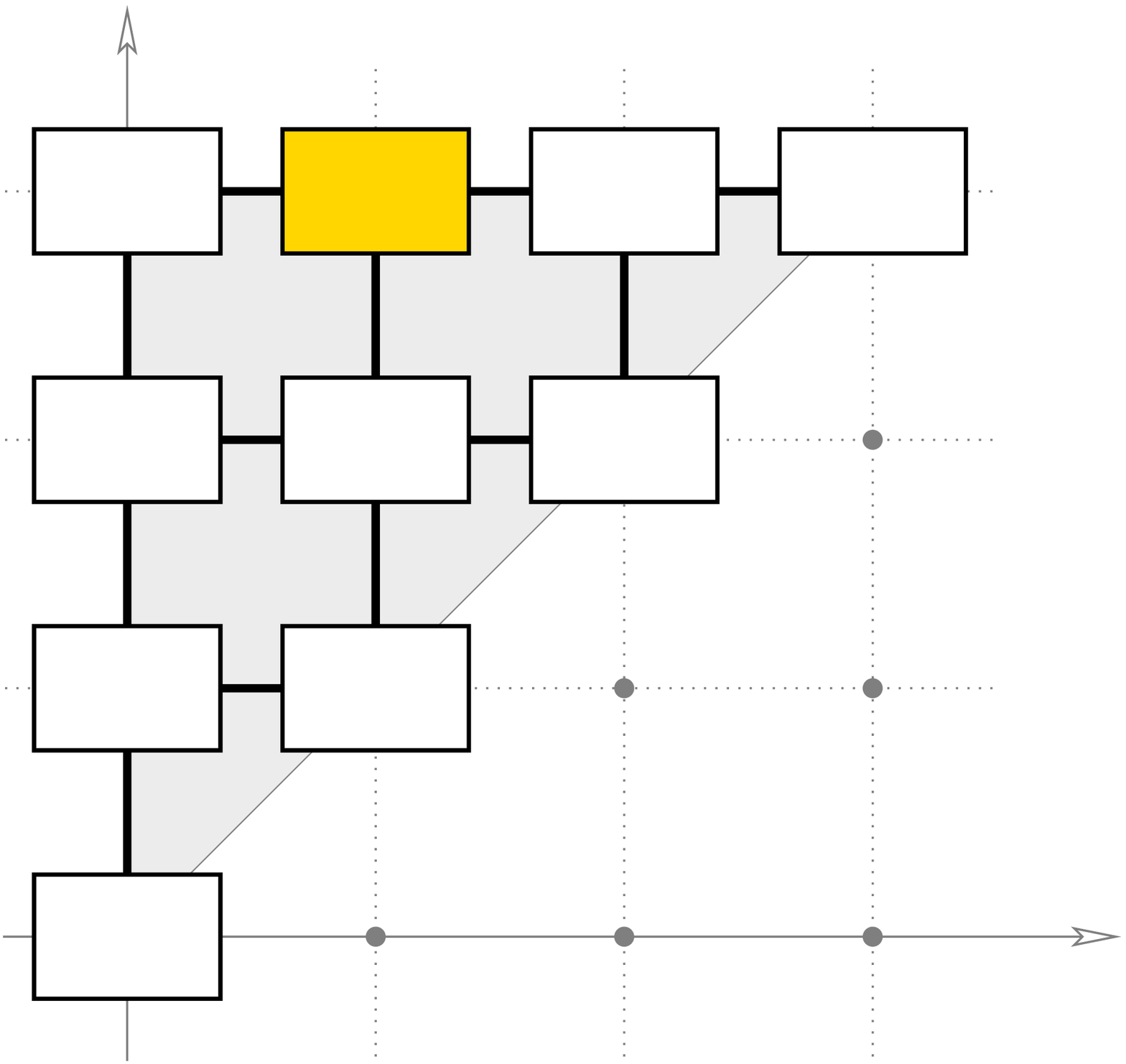}
    \put(4,78){\small$10001$}  \put(26,78){\small$01001$}  \put(48,78){\small$00101$}  \put(70,78){\small$00011$}
    \put(8,73){\small$0,\!3$}  \put(30,73){\small$1,\!3$}  \put(52,73){\small$2,\!3$}  \put(74,73){\small$3,\!3$}
    \put(4,56){\small$10010$}  \put(26,56){\small$01010$}  \put(48,56){\small$00110$}
    \put(8,51){\small$0,\!2$}  \put(30,51){\small$1,\!2$}  \put(52,51){\small$2,\!2$}
    \put(4,34){\small$10100$}  \put(26,34){\small$01100$}
    \put(8,29){\small$0,\!1$}  \put(30,29){\small$1,\!1$}
    \put(4,12){\small$11000$}
    \put(8,7){\small$0,\!0$}
  \end{overpic}
  \caption{The facet $01001$ of $\stc(\Delta_2\times\Delta_3)$ and the dual graph of
    $\stc(\Delta_2\times\Delta_3)$ with the facet $01001$ marked. \label{fig:staircase}}
\end{figure}

We are interested in one particular triangulation of $\Delta_m\times\Delta_n$, the \emph{staircase
  triangulation} $\stc_{m,n}=\stc(\Delta_m\times\Delta_n)$, which can be described as follows.
Consider a rectangular grid of size $m+1$ by $n+1$.  Each node in the grid corresponds to one vertex
of $\Delta_m\times\Delta_n$.  The facets of $\stc_{m,n}$, described as subsets of these nodes,
correspond to the non-descending and not-returning paths from the lower left node to the upper right
node.  These paths, which go only right or up, but never left nor down, look like staircases, and
hence the name; see \mbox{Figure~\ref{fig:staircase} (left)}.

The choice of ``right'' and ``up'' in the definition of $\stc_{m,n}$ implicitly assumes an ordering
of the vertices of both factors.  Throughout this paper we will keep this ordering fixed.  The
staircase triangulation of $\Delta_m\times\Delta_n$ is the same as the placing triangulation induced
by the \emph{lexrev ordering}, that is, the lexicographic ordering of the vertices with the reversed
ordering of the vertices of the second factor.  In particular, $\stc_{m,n}$ is a regular
triangulation.

Each such staircase can be encoded as a \emph{shuffle} of ``up'' and ``right'' moves.  The name
``shuffle'' reflects the fact that the number of ``up'' and ``right'' moves is always the same, but
their order is all that matters.  We write the shuffle in Figure~\ref{fig:staircase} as the
bit-string $01001$, where $0$ means ``up'' and $1$ means ``right''.  The staircase triangulations
occurred in Eilenberg and Steenrod~\cite[Section~II.8]{MR0050886}; see also Billera, Cushman, and
Sanders~\cite{MR976522}, Gel$'$fand, Kapranov, and \mbox{Zelevinsky~\cite[\S7.D]{MR1264417}}, and
Santos~\cite{MR1758756}.

Yet another way to encode a facet~$F$ of~$\stc_{m,n}$ is to assign a vector $s(F) \in \NN^m$ as
follows.  The bit-string $11\dots 100\dots0$ corresponds to the origin, and for an arbitrary facet
$F$ the $k$-th entry $s(F)_k$ measures the difference between the position of the $k$-th one in the
bit-representation of~$F$ and $k$.  This difference may be viewed as the number of ``shifts to the
right'' of the $k$-th one, starting with the bit-string corresponding to the origin.  For example,
the bit-string~$01001$ in Figure~\ref{fig:staircase} is mapped to~$(1,3)$.

Via the map $s$ the facets of $\stc_{m,n}$ correspond to the integer points in the polytope
\[\mathcal{S}_{m,n} = \SetOf{s\in\RR^m}{
  \begin{aligned}
    0 \leq s_k \leq n \; & \text{for $1 \leq k \leq m$} \\
    s_k \leq s_l \; & \text{for $k<l$}
  \end{aligned}} \; .\]

This provides us with a convenient description of the dual graph of~$\stc_{m,n}$.
Let~$\mathcal{L}_m$ be the $m$-dimensional cubic grid, that is, the infinite graph with node set
$\ZZ^m$, and two nodes are adjacent if they differ in exactly one coordinate by one.

We denote the dual graph of a simplicial complex~$K$ by $\Gamma^*(K)$.  Its nodes are the facets
of~$K$ and two facets are adjacent if they differ in one vertex.

\begin{prop} \label{prop:graph_as_intersection}
  The dual graph $\Gamma^*(\stc_{m,n})$ is the subgraph of~$\mathcal{L}_m$ induced by the node set
  $\mathcal{S}_{m,n}\cap\ZZ^m$.  In particular, this graph is bipartite.
\end{prop}

To conclude this section we mention further aspects of the staircase triangulations, which are,
however, inessential for the understanding of rest of this paper.

\begin{rem}
  Bit-strings of length~$m+n$ with precisely $m$ ones correspond to the vertices of the hypersimplex
  $H(m+n,m)$.  The graph $\Gamma^*(\stc_{m,n})$ is a (not induced) subgraph of the vertex-edge graph
  of~$H(m+n,m)$.
  The Cayley trick establishes a one-to-one correspondence between the regular triangulations of
  $\Delta_m\times\Delta_n$ and the fine mixed subdivisions of $(n+1)\Delta_m$; see
  Santos~\cite{MR2134766}.
  In a different context regular triangulations of $\Delta_m\times\Delta_n$ recently re-appeared as
  the tropical convex hulls of $n+1$ points in the tropical projective space $\tropPG^m$; see
  Develin and Sturmfels~\cite{TropConvex}.  The staircase triangulations arise as the tropical
  cyclic polytopes of Block and Yu~\cite{math.MG/0503279}.
\end{rem}

\section{Products of Simplicial Complexes}

Let~$K$ and~$L$ be two abstract simplicial complexes. Then the product space $\abs{K}\times\abs{L}$
is equipped with the structure of a cell complex whose cells are the products $f\times g$, where $f$
is a face of~$K$ and $g$ is a face of~$L$.  This section is about the study of triangulations
of~$\abs{K}\times\abs{L}$ which refine this natural cell structure.

\subsection{The Simplicial Product}

Assume that $\dim K=m$ and $\dim L=n$, and denote the vertex sets of $K$ and $L$ by $V_K$ and $V_L$,
respectively.  We choose a linear ordering $O_K$ of $V_K$ and another linear ordering $O_L$ of
$V_L$.  The product $O_K\times O_L$, defined by
\[
(v,w)\geq (v',w')\;\Leftrightarrow\;
v\geq v'\;\text{and}\;  w\geq w'\;,
\]
is a partial ordering of the set $V_K\times V_L$.  Let $\pi_K:V_K \times V_L \rightarrow V_K$ and
$\pi_L:V_K \times V_L \rightarrow V_L$ be the canonical projections.

We define the \emph{simplicial product} (with respect to the vertex orderings~$O_K$ and~$O_L$)
of~$K$ and~$L$ as
\[ K \simpp L = \SetOf{F\subseteq  V_K\times V_L}{
  \begin{matrix}
    \pi_K(F) \in K \; \text{and} \; \pi_L(F) \in L\,,\\
    \text{and} \; O\left|_F\right. \; \text{is a total ordering}
  \end{matrix}} \; .\]

The simplicial product $K \simpp L$ appeared earlier in Eilenberg and
Steenrod~\cite[Section~II.8]{MR0050886} as the ``Cartesian product'', and in
Santos~\cite{MR1758756}, who calls it the ``staircase refinement''. Both sources prove the
staircase triangulation to be a triangulation of the space $\abs{K}\times\abs{L}$ on the vertex set
$V_K\times V_L$.

\begin{figure}[htbp]\centering
  \includegraphics[width=.9\textwidth]{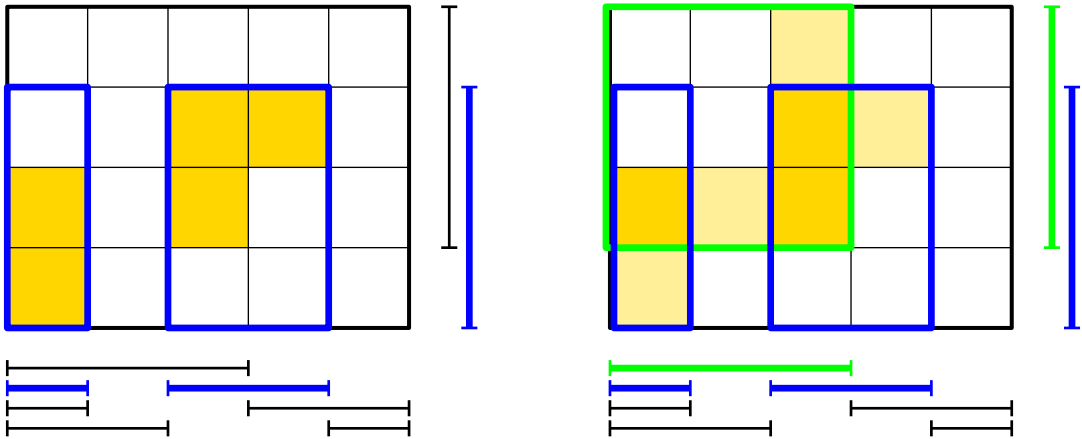}
  \caption{A facet defining path of the simplicial product of two different triangulations of the
    square. On the right two facets intersecting in a low dimensional
    face.\label{fig:product_triang}}
\end{figure}

Let $k = |V_K|$ and $l = |V_L|$ denote the number of vertices of~$K$ and~$L$, respectively.  There
is a convenient way to visualize the simplicial product in the $(k\times l)$-grid~$\cR$: We label
the columns of~$\cR$ with the vertices of~$K$ according to the vertex order~$O_K$, and we label the
rows of~$\cR$ with the vertices of~$L$ according to the vertex order~$O_L$.  For each~$f \in K$
and~$g \in L$ let~$\cR_{f,g}$ be the minor of~$\cR$ induced by~$f$ and~$g$. Then we may think of the
facets of the simplicial product as the collection of all ascending paths in~$\cR_{f,g}$ starting
bottom-left and finishing top-right.  This is a direct generalization of the staircase triangulation
of the product of two simplices; see Figure~\ref{fig:product_triang}.  More precisely, we may view
the simplicial product $K\simpp L$ as a subcomplex of the staircase triangulation of the product of
a $(k-1)$-simplex and an $(l-1)$-simplex.

\begin{figure}[htbp]\centering
  \begin{overpic}[width=.95\textwidth]{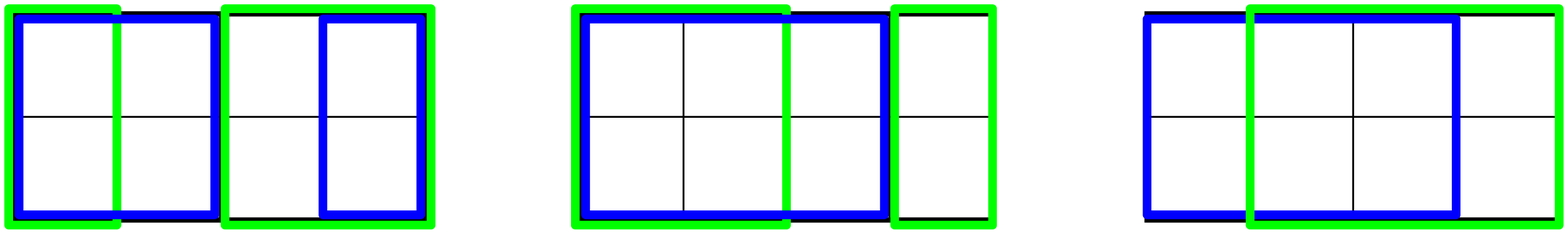}
    \put(29,11){\small $a$}
    \put(29,17){\small $b$}
    \put(3,5){\small $1$}
    \put(10,5){\small $0$}
    \put(17,5){\small $3$}
    \put(23,5){\small $2$}
    \put(65,11){\small $a$}
    \put(65,17){\small $b$}
    \put(39,5){\small $1$}
    \put(46,5){\small $2$}
    \put(53,5){\small $0$}
    \put(59,5){\small $3$}
    \put(101,11){\small $a$}
    \put(101,17){\small $b$}
    \put(76,5){\small $0$}
    \put(82,5){\small $1$}
    \put(89,5){\small $2$}
    \put(96,5){\small $3$}
  \end{overpic}

  \mbox{
    \begin{overpic}[width=.29\textwidth]{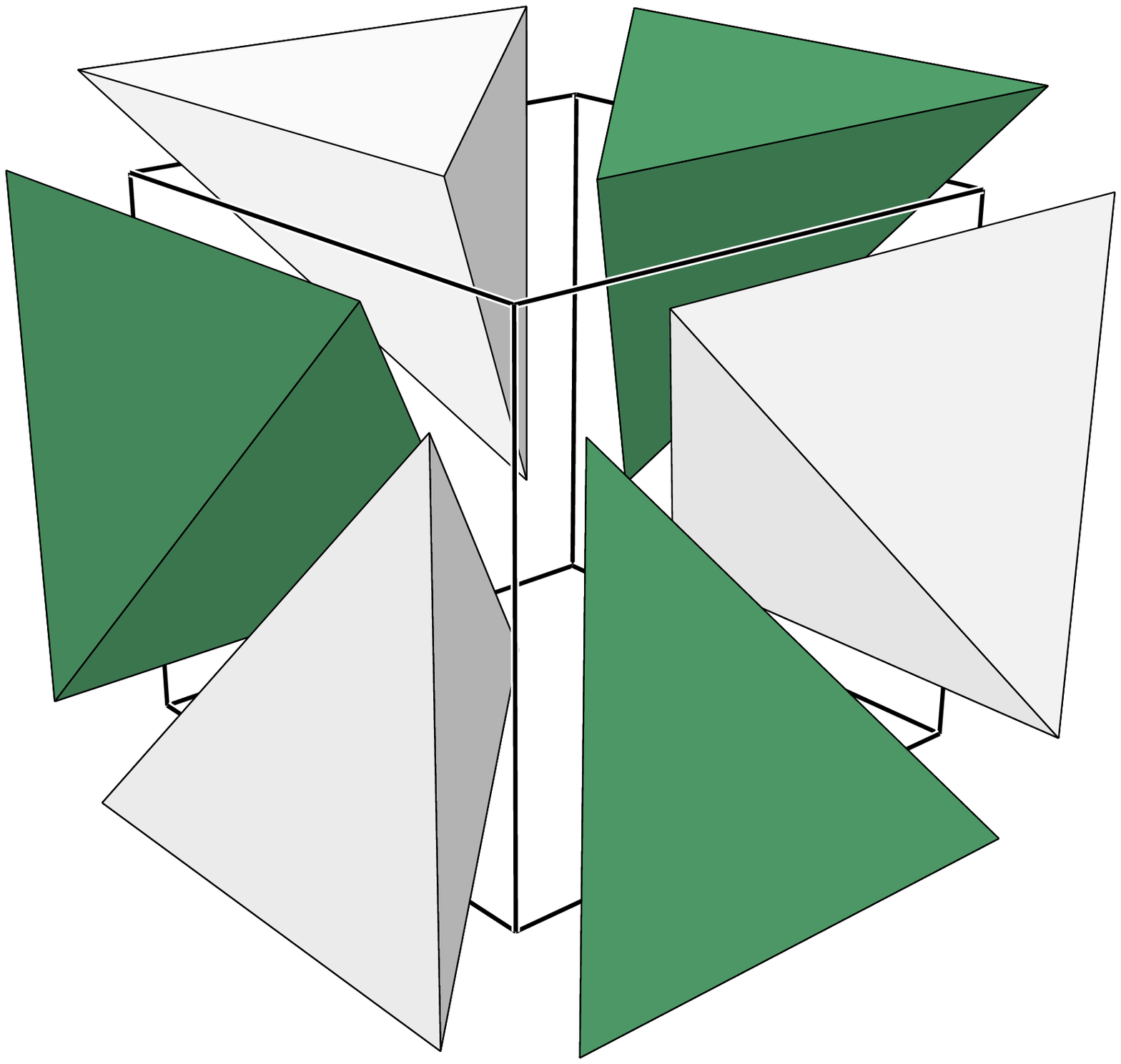}
      \put(-6,72){\small $(3,a)$}
      \put(44,84){\small $(2,a)$}
      \put(91,71){\small $(0,a)$}
      \put(-4,21){\small $(3,b)$}
      \put(39,-6){\small $(1,b)$}
      \put(87,18){\small $(0,b)$}
    \end{overpic}
    
    \hspace{0.6cm}
    \begin{overpic}[width=.29\textwidth]{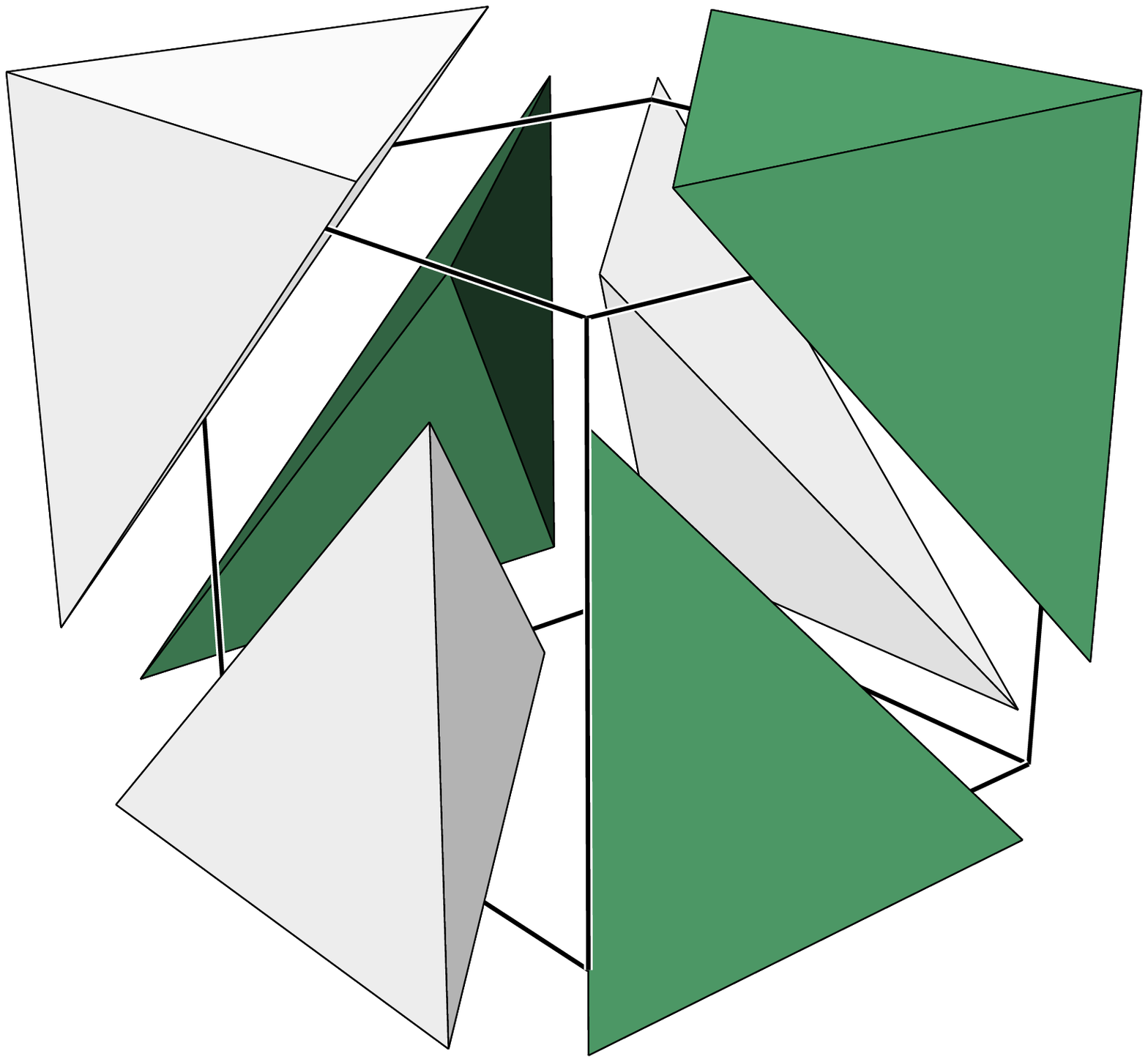}
    \end{overpic}
    
    \hspace{0.2cm}
    \begin{overpic}[width=.29\textwidth]{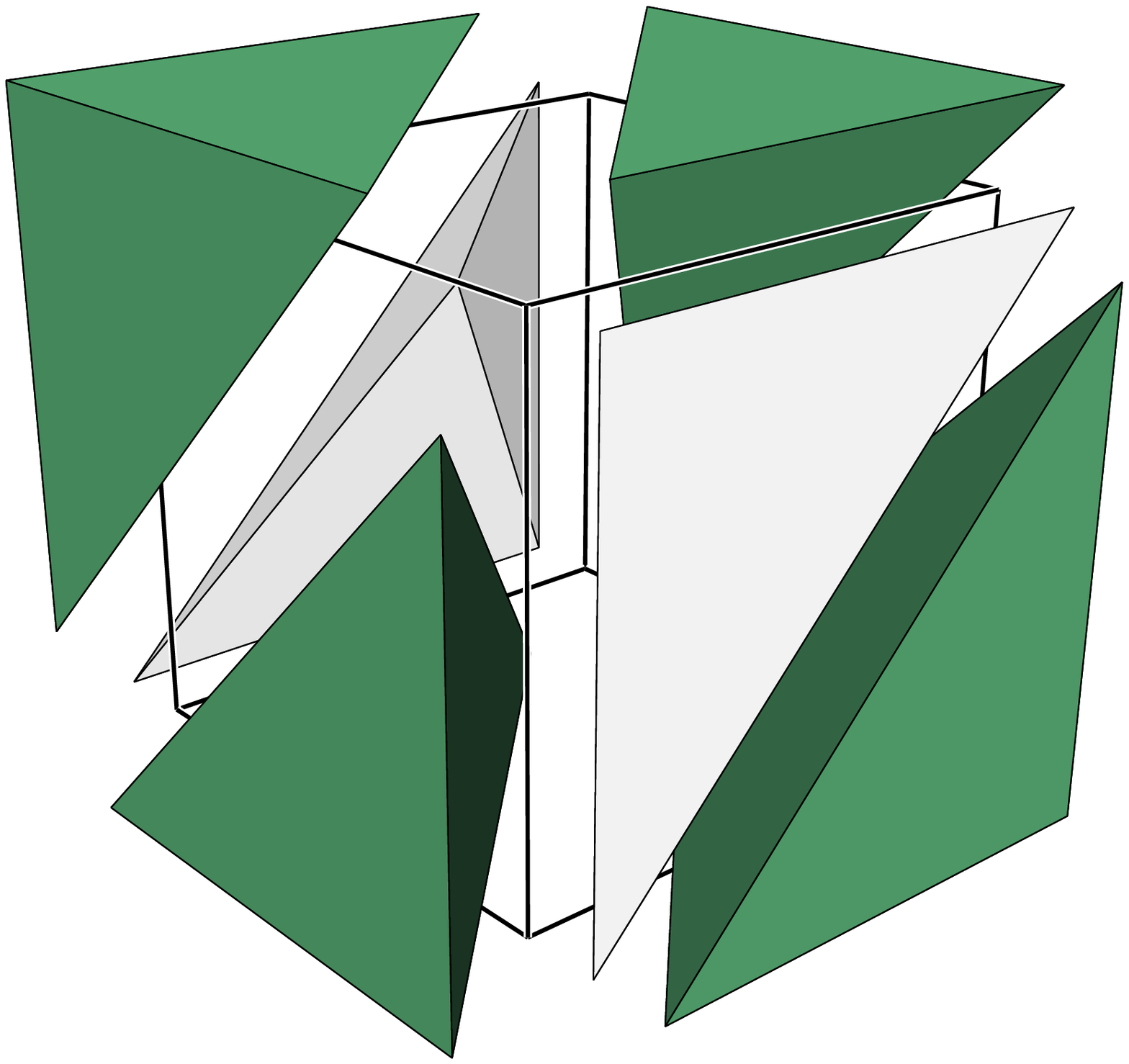}
    \end{overpic}
  }
  
  \includegraphics[width=.01cm]{empty_box.eps}
  \caption{Three different orderings of the vertices of the triangulated square
    $\{\{0,1,2\},\{1,2,3\}\}$ and the resulting
    regular triangulations of the $3$-cube.  The vertices $0$ and $3$ of the square are colored the
    same, and the top-front vertex of the $3$-cube is labeled $(1,a)$, and the bottom-back vertex is
    labeled $(2,b)$.  The second and third $3$-cube are labeled the same. \label{fig:three_orderings}}
\end{figure}

The ordering of the vertices of~$K$ and~$L$ is crucial to $K \simpp L$.
Figure~\ref{fig:three_orderings} depicts the product of the triangulated unit square with the unit
interval. The three distinct orderings of the vertices of the triangulated square yield three
pairwise non-isomorphic triangulations of the $3$-cube $C_3$ decomposed as $C_2 \times I$.

\subsection{Foldable Simplicial Complexes} \label{subsec:foldable_complexes}

An $m$-dimensional pure simplicial complex $K$ is called \emph{foldable} if $K$ admits a
non-degenerate simplicial map to an $m$-simplex.  Equivalently, the $1$-skeleton of~$K$ is
$(m+1)$-colorable in the graph-theoretic sense: that is, there is a map $c$ from the vertex set $V$
to the set $[m+1]$ such that for each $1$-face $\{u,v\}\in K$ we have $c(u)\ne c(v)$.  Here
$[k]=\{0,\ldots ,k-1\}$ denotes the set of the first $k$ integers.  Notice that there is no coloring
of the vertices of~$K$ with less than $m+1$ colors, since the $m+1$ vertices of any facet form a
clique.  If $K$ is foldable with a connected dual graph then the $(m+1)$-coloring of~$K$ is unique
up to renaming the colors.

Goodman and Onishi~\cite{MR515556} observed that the $4$-Color-Theorem is equivalent to the property
that each simplicial $3$-polytope admits a foldable triangulation (with or without additional
vertices in the interior). 

\begin{rem}
  Other sources, including Billera and Bj\"orner~\cite{billera_bjoerner}, Stanley~\cite{MR725505},
  Soprunova and Sottile~\cite{math.AG/0409504}, and \cite{MR1900311,MR1967999}, call foldable
  simplicial complexes ``balanced.''  However, this seems to create conflicts with other concepts: A
  triangulation of a polygon whose dual graph is a balanced tree is sometimes called ``balanced'',
  and a minimal set of affinely dependent vertices of a polytope with an equal number of positive
  and negative coefficients is called a ``balanced'' circuit in Bayer~\cite{MR1231196}. Goodman and
  Onishi call foldable triangulations (of balls and spheres) ``even.''  However, this does not
  describe the situation in the non-simply connected case.  For these reasons we suggest the name
  ``foldable'' instead.
\end{rem}

If $K$ is pure and, additionally, certain global and local connectivity assumptions are satisfied,
then $K$ is foldable if and only if its group of projectivities is trivial.  These connectivity
assumptions hold, for instance, when $K$ is the triangulation of a manifold (with or without
boundary).  Moreover, in this case, foldability implies that the dual graph of $K$ is bipartite.
The converse holds for simply connected combinatorial manifolds.  For these facts and related
results see~\cite{MR1900311,MR1967999}.  In the following we study products of foldable simplicial
complexes.

Let $[k]=\{0,\dots,k-1\}$ be the vertex set of $K$.  Assume that there is a coloring of $K$ given by
a weakly monotone map $c_K:[k]\to[m+1]$.  Then we call the natural ordering on~$[k]$ \emph{color
  consecutive}.  Any foldable complex admits (many) color consecutive orderings.

\begin{prop} \label{prop:foldable_product}
  If~$K$ and~$L$ are foldable simplicial complexes with color consecutive vertex orderings then the
  corresponding simplicial product $K\simpp L$ is foldable.
\end{prop}

\begin{proof}
  Let the vertex sets of $K$ and $L$ be $[k]$ and $[l]$, respectively, with weakly monotone coloring
  maps $c_K:[k]\to[m+1]$ and $c_L:[l]\to[n+1]$. We define
  \[
  c:[k]\times[l]\to[m+n+1]:(v,w)\mapsto c_K(v)+c_L(w)\;.
  \]
  In order to show that $c$ is a coloring of $K\simpp L$ it suffices to check that
  each facet contains each color at most once.  Each facet~$F$ of $K\simpp L$ is contained in a
  unique cell $f\times g$ where $f$ is a facet of $K$ and $g$ is a facet of $L$. Let $v\times w$ and
  $v' \times w'$ be distinct vertices of~$F$. We may assume $v < v'$; then $w \leq w'$ since~$F$ is
  a facet of the staircase triangulation of $f\times g$. As the restrictions $c_K \: \restr_f$
  and~$c_L \: \restr_g$ are strictly monotone we have
  $c(v,w)=c_K(v)+c_L(w)<c_K(v')+c_L(w')=c(v',w')$.  For an example see
  Figure~\ref{fig:lifted_product}.
\end{proof}

In what follows below it is essential that it is not necessary to have color consecutive orderings
for the factors in order to obtain a foldable simplicial product triangulation.

\begin{exmp}\label{exmp:bipyramid}
  Let $B_n$ be the triangulation of the bipyramid over the $(n-1)$-simplex~$\Delta_{n-1}$ formed of
  two $n$-simplices sharing a facet.  Combinatorially, $B_n$ is the join of $\Delta_{n-1}$ with the
  zero-dimensional sphere $\Sph^0$ consisting of two isolated points.  The triangulation $B_n$ is
  obviously foldable.  The \emph{symmetric vertex ordering} $S_n$ on $B_n$ starts with one of the
  two apices and ends with the other apex, the vertices of $\Delta_{n-1}$ come in between.  That is
  to say, we take $[n+2]$ as the vertex set of~$B_n$, where $0$ and $n+1$ are the apices, and a
  coloring map $s_n:[n+2]\to[n+1]:w\mapsto w\mod (n+1)$.  Because of the symmetry properties of
  $B_n$ the precise ordering of the vertices $1,2,\dots,n$ does not matter.  Likewise it is not
  necessary to distinguish the two apices.
\end{exmp}

The triangulation $B_n$ with the symmetric vertex ordering will be used in the construction of
certain cube triangulations in Section~\ref{sec:cubes}.

\begin{prop} \label{prop:symmetric_ordering}
  Let~$K$ be a foldable simplicial complex with a color consecutive ordering~$O_K$.  Then the
  simplicial product $K\simpp B_n$ with respect to $O_K$ and~$S_n$ is foldable.
\end{prop}

\begin{proof}
  We use almost the same coloring scheme as in Proposition~\ref{prop:foldable_product}.  Let~$[k]$
  be the vertex set of $K$, and let $c_K:[k]\to[m+1]$ be a weakly monotone coloring map.  We define
  \[c:[k]\times[n+2]\to[m+n+1]:(v,w)\mapsto c_K(v)+w\mod (m+n+1).\]  This, indeed, is a coloring since
  there is no facet of $K\simpp B_n$ containing both, a vertex of the type $(v,0)$ and a vertex of
  the type $(v,n+1)$.
\end{proof}

We refer to Figure~\ref{fig:three_orderings} for the three different simplicial products of an interval
with a square arising from the two color consecutive and the symmetric vertex ordering of the square
(which is a bipyramid over a $1$-simplex).

\subsection{Regular Triangulations of Polytopes}

Let $P$ be an $m$-dimensional convex polytope in~$\RR^m$, and let $K$ be a triangulation of~$P$ with
vertex set~$V$.  The triangulation $K$ is \emph{regular} if there is a convex function
$\lambda:\RR^m\to\RR$ such that $K$ coincides with the polyhedral subdivision of~$P$ induced by the
lower convex hull of the set $\SetOf{(v,\lambda(v))\in\RR^{m+1}}{v\in V}$.  In this case $\lambda$
is called a \emph{lifting function} for~$K$.  Since we want to stress that a regular triangulation
only depends on~$P$ and~$\lambda$ we denote such a triangulation as~$P^\lambda$.

Choose (pairwise distinct) points $p_1,\dots,p_k$ in~$P$ such that $\conv\{p_1,\dots,p_k\}=P$.
This implies that the vertices of $P$ occur among the chosen points.  Then the \emph{placing
  triangulation} of~$P$ with respect to the chosen points in the given ordering is the regular
triangulation of $P$ with vertex set $\{p_1,\dots,p_k\}$ and a lifting function $\lambda$ such that
$(p_l,\lambda(p_l))$ is above all affine hyperplanes spanned by points in the set
$\{(p_1,\lambda(p_1)),\dots,(p_{l-1},\lambda(p_{l-1}))\}$.  A point $(p,\lambda(p))$ 
lies \emph{above} the affine hyperplane $H\subset\RR^{m+1}$ spanned by the points
$\{(p_1,\lambda(p_1)),\dots,(p_{m+1},\lambda(p_{m+1}))\}$ if and only if the unique $\lambda' \in
\RR$ with
\begin{equation}\label{eq:lies_above}
\det
\begin{pmatrix}
  1 & 1 & 1 & \dots & 1\\
  p & p_1 & p_2 & \dots & p_{m+1}\\
  \lambda' & \lambda(p_1) & \lambda(p_2) & \dots & \lambda(p_{m+1})
\end{pmatrix}
= 0
\end{equation}
satisfies $\lambda' < \lambda(p)$.

\begin{exmp}\label{exmp:lifting_function_stc}
  Consider the standard simplices $\Delta_m = \conv\{0,e_1,\dots,e_m\}$ and $\Delta_n =
  \conv\{0,e_1,\dots,e_n\}$.  To simplify the formulae below we set $e_0=0$.  Then the lexrev
  ordering on the vertices of the product $\Delta_m\times\Delta_n$ is given as
  \[
  O: \{e_0,\dots,e_m\} \times \{e_0,\dots,e_n\} \to [(m+1)(n+1)]:
  (e_i,e_j) \mapsto (n+1)i+(n-j) \; .
  \]
  Applying Equation~\ref{eq:lies_above} and an easy computation shows that
  \[
  \omega: \{e_0,\dots,e_m\} \times \{e_0,\dots,e_n\} \to \RR: (v,w) \mapsto 2^{O(v,w)}
  \]
  is a lifting function for the staircase triangulation, that is,
  $(\Delta_m\times\Delta_n)^\omega=\stc_{m,n}$.  Additionally, this shows that $\stc_{m,n}$ is a
  placing triangulation with respect to the lexrev ordering.
\end{exmp}

\begin{prop}\label{prop:regular}
  Let $P^\lambda$ and $Q^\mu$ be regular triangulations of an $m$-polytope \mbox{$P\subset\RR^m$} and an
  $n$-polytope $Q\subset\RR^n$, respectively.  Then the simplicial product $P^\lambda\simpp Q^\mu$
  is a regular triangulation of the polytope $P\times Q$ for any vertex orderings~$O_{P^\lambda}$
  and~$O_{Q^\mu}$.
\end{prop}

\begin{figure}[t]\centering
  \psfrag{0}[][]{\hspace{-0.2 cm} $\mathbf{0}$}
  \psfrag{1}[][]{\hspace{-0.2 cm} $\mathbf{4}$}
  \psfrag{2}[][]{\hspace{-0.2 cm} $\mathbf{1}$}
  \psfrag{3}[][]{\hspace{-0.2 cm} $\mathbf{8}$}
  \psfrag{4}[][]{\hspace{-0.2 cm} $\mathbf{5}$}
  \psfrag{5}[][]{\hspace{-0.2 cm} $\mathbf{2}$}
  \psfrag{6}[][]{\hspace{-0.2 cm} $\mathbf{12}$}
  \psfrag{7}[][]{\hspace{-0.2 cm} $\mathbf{9}$}
  \psfrag{8}[][]{\hspace{-0.2 cm} $\mathbf{6}$}
  \psfrag{9}[][]{\hspace{-0.2 cm} $\mathbf{3}$}
  \psfrag{10}[][]{\hspace{-0.2 cm} $\mathbf{13}$}
  \psfrag{11}[][]{\hspace{-0.2 cm} $\mathbf{10}$}
  \psfrag{12}[][]{\hspace{-0.2 cm} $\mathbf{7}$}
  \psfrag{13}[][]{\hspace{-0.2 cm} $\mathbf{14}$}
  \psfrag{14}[][]{\hspace{-0.2 cm} $\mathbf{11}$}
  \psfrag{15}[][]{\hspace{-0.2 cm} $\mathbf{15}$}
  \psfrag{a}[][]{\hspace{-0.2 cm} $\mathbf{b}$}
  \psfrag{b}[][]{\hspace{-0.2 cm} $\mathbf{c}$}
  \psfrag{c}[][]{\hspace{-0.2 cm} $\mathbf{a}$}
  \psfrag{d}[][]{\hspace{-0.2 cm} $\mathbf{d}$}
  \psfrag{aa}[][]{\hspace{-0.2 cm} $\mathbf{d'}$}
  \psfrag{bb}[][]{\hspace{-0.2 cm} $\mathbf{b'}$}
  \psfrag{cc}[][]{\hspace{-0.2 cm} $\mathbf{c'}$}
  \psfrag{dd}[][]{\hspace{-0.2 cm} $\mathbf{a'}$}
  \includegraphics[width=.85\textwidth]{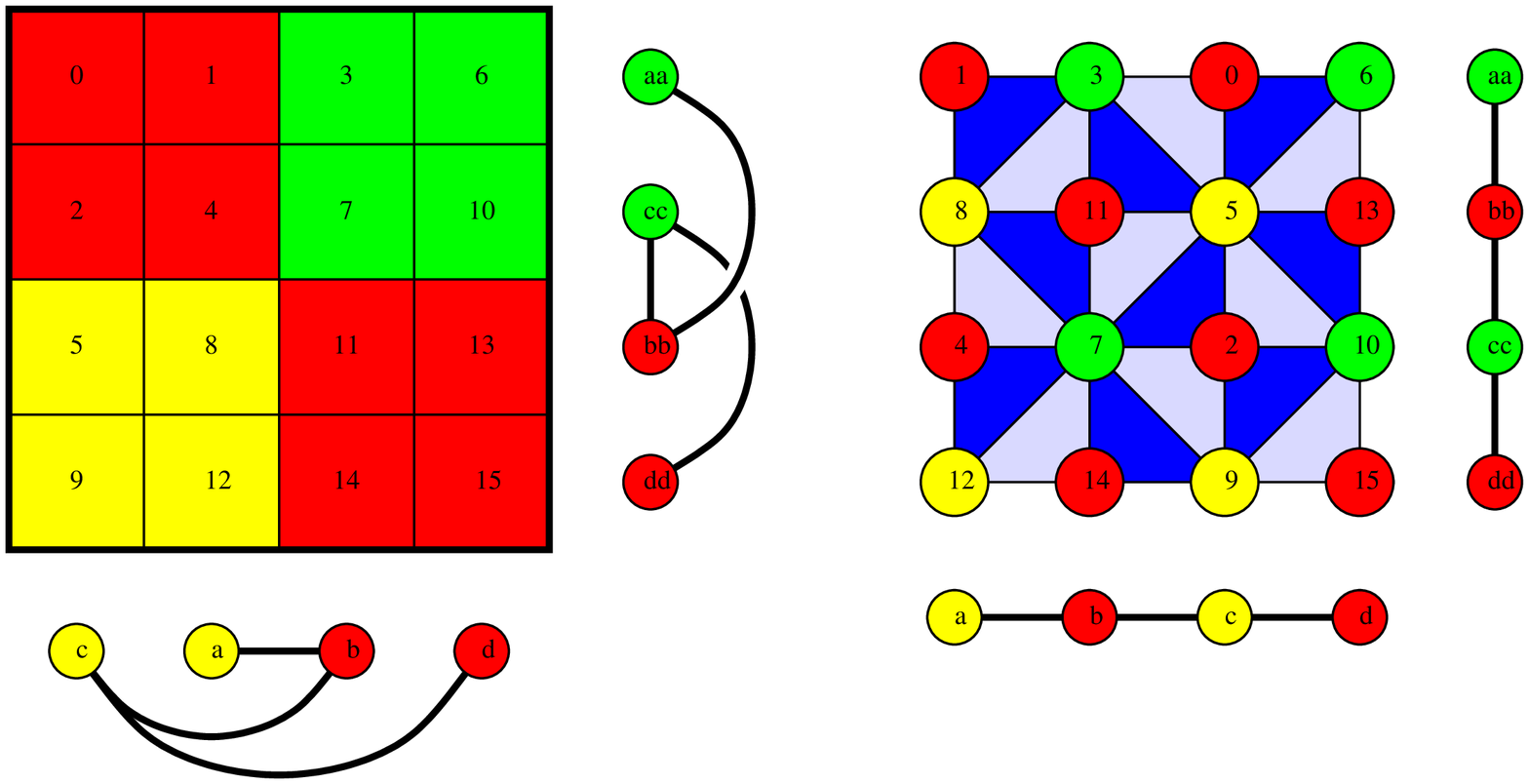}

  \psfrag{0}[][]{\hspace{-0.6 cm} $\mathbf{4}$}
  \psfrag{1}[][]{\hspace{-0.6 cm} $\mathbf{0}$}
  \psfrag{2}[][]{\hspace{-0.6 cm} $\mathbf{5}$}
  \psfrag{3}[][]{\hspace{-0.6 cm} $\mathbf{8}$}
  \psfrag{4}[][]{\hspace{-0.6 cm} $\mathbf{1}$}
  \psfrag{5}[][]{\hspace{-0.6 cm} $\mathbf{6}$}
  \psfrag{6}[][]{\hspace{-0.68 cm} $\mathbf{12}$}
  \psfrag{7}[][]{\hspace{-0.6 cm} $\mathbf{9}$}
  \psfrag{8}[][]{\hspace{-0.6 cm} $\mathbf{2}$}
  \psfrag{9}[][]{\hspace{-0.6 cm} $\mathbf{7}$}
  \psfrag{10}[][]{\hspace{-0.68 cm} $\mathbf{13}$}
  \psfrag{11}[][]{\hspace{-0.68 cm} $\mathbf{10}$}
  \psfrag{12}[][]{\hspace{-0.6 cm} $\mathbf{3}$}
  \psfrag{13}[][]{\hspace{-0.68 cm} $\mathbf{14}$}
  \psfrag{14}[][]{\hspace{-0.68 cm} $\mathbf{11}$}
  \psfrag{15}[][]{\hspace{-0.68 cm} $\mathbf{15}$}
  \includegraphics[width=.65\textwidth]{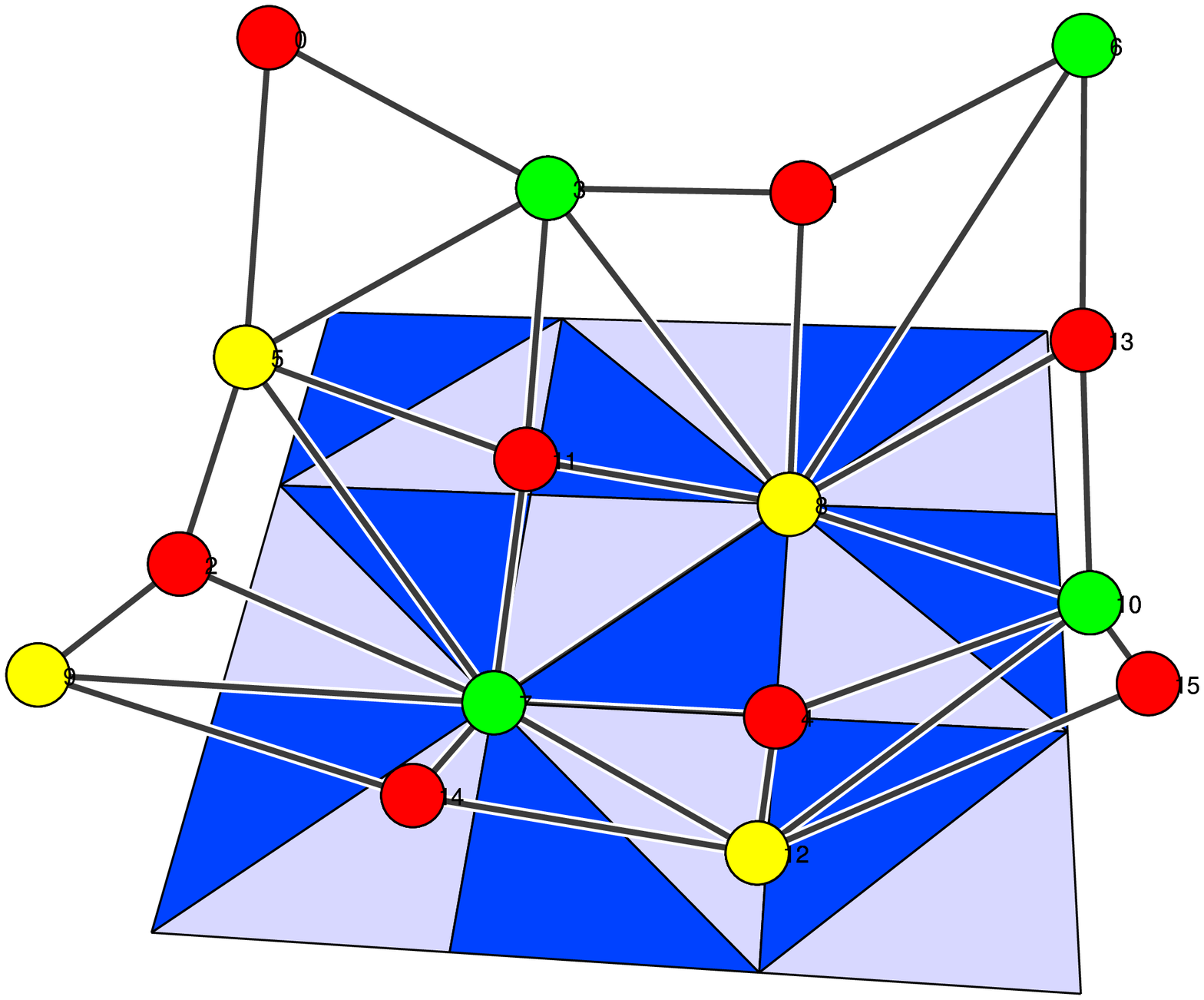}
  \caption{Simplicial product of a path $I$ of length $3$ with itself, using color consecutive
    vertex orderings. The vertices of the product are colored according to the color scheme from the
    proof of Proposition~\ref{prop:foldable_product} and are labeled in lexrev order.}
  \label{fig:lifted_product}
\end{figure}

\begin{proof}
  Let $V_{P^\lambda}$ be the vertex set of~$P^\lambda$ equipped with a linear ordering
  $O_{P^\lambda}$, and let~$V_{Q^\mu}$ be the vertex set of~$Q^\mu$ with a linear ordering
  $O_{Q^\mu}$.  The simplicial product $P^\lambda\simpp
  Q^\mu$ (with respect to $O_{P^\lambda}$ and $O_{Q^\mu}$) is a triangulation of the product
  $P\times Q$ on the vertex set~$V_{P^\lambda}\times V_{Q^\mu}$.
  
  Let $\lambda: V_{P^\lambda}\to \RR$ and $\mu: V_{Q^\mu} \to \RR$ be lifting
  functions of~$P^\lambda$ and~$Q^\mu$.  We construct a lifting function $\omega:
  V_{P^\lambda}\times V_{Q^\mu} \to \RR$ of $P^\lambda \simpp Q^\mu$ in two steps.  First
  consider the map
  \[
  \omega_0:V_{P^\lambda}\times V_{Q^\mu}\to\RR:(x,y)\mapsto \lambda(x)+\mu(y)\; , 
  \]
  which is a
  lifting function for the polytopal complex $P^\lambda\times Q^\mu$.  In the second step~$\omega_0$
  has to be perturbed such that it induces a staircase triangulation on each cell of $P^\lambda
  \times Q^\mu$.  To this end recall that the staircase triangulations are placing, and that the lexrev
  ordering~$O$ on $V_{P^\lambda}\times V_{Q^\mu}$ induces a placing order on each product of
  simplices $f\times g$ where $f\in P^\lambda$ and $g\in Q^\mu$.  Now define~$\omega$ as an
  $\epsilon$-perturbation of $\omega_0$ by the lifting function from
  Example~\ref{exmp:lifting_function_stc} corresponding to~$O$:
  \begin{equation} \label{eq:omega1}
    \omega:  V_{P^\lambda}\times V_{Q^\mu} \to \RR : (v,w) \mapsto \lambda(v) + \mu(w)
    + \epsilon 2^{O(v,w)} \; ,
  \end{equation}
  for a sufficiently small $\epsilon>0$.  Viewing the simplicial product again as subcomplex of the
  staircase triangulation of two large simplices, shows that, indeed $(P\times
  Q)^\omega=P^\lambda\simpp Q^\mu$.  For an example see Figure~\ref{fig:lifted_product}.
\end{proof}

In general, there may be several perturbations which lead to different lifting functions but which
induce the same triangulations.  An important special case occurs if the triangulations $P^\lambda$
and $Q^\lambda$ additionally are foldable.  In this case it is possible to define a perturbation
which only depends on the color classes of the vertices of the factors:

\begin{exmp}
  Let $c_{P^\lambda} : V_{P^{\lambda}} \to [m+1]$ and $c_{Q^\mu} : V_{Q^{\mu}} \to [n+1]$ be
  coloring maps.  Using color consecutive vertex orderings for $V_{P^{\lambda}}$ and $V_{Q^{\mu}}$
  and the resulting lexrev ordering~$O$ for the vertices of $P \simpp Q$ we may choose a different
  perturbation than in Equation~\eqref{eq:omega1}.  This yields the following lifting function
  \begin{equation} \label{eq:omega2}
    \omega:  V_{P^\lambda}\times V_{Q^\mu} \to \RR : (v,w) \mapsto \lambda(v) + \nu(w)
    + \epsilon \: 2^{(n+1)c_{P^\lambda}(v)+(n-c_{Q^\mu}(w))} \; ,
  \end{equation}
  for $\epsilon>0$ sufficiently small.  Note that we use the same perturbation $\epsilon 2^{(n+1)i+(n-j)}$ for all
  vertices $(v,w)$ with $c_{P^\lambda}(v)=i$ and $c_{Q^\mu}(w)=j$.  Let us restrict our attention to
  a cell $f\times g$ for facets $f\in P^{\lambda}$ and $g\in Q^\mu$.  Since any color $i \in [m+1]$
  appears exactly once in the coloring of~$f$ and any color $j \in [n+1]$ appears exactly once in
  the coloring of~$g$, respectively, there is exactly one vertex $(v,w) \in f\times g$ with
  $c_{P^\lambda}(v)=i$ and $c_{Q^\mu}(w)=j$ for each $(i,j)\in [m+1]\times [n+1]$.  Hence $\omega$
  restricted to $f\times g$ induces the staircase triangulation $f\simpp g$ from
  Example~\ref{exmp:lifting_function_stc}, and $\omega$ induces the simplicial product triangulation
  $(P\times Q)^\omega=P^\lambda\simpp Q^\mu$ on $P^\lambda\times Q^\mu$.
\end{exmp}

\section{Triangulations of Lattice Polytopes}

Let $P$ be an $m$-dimensional \emph{lattice polytope}, that is, we assume that its vertex
coordinates are integral. Since the determinant of an integral matrix is an integer it follows that
the \emph{normalized volume} $\nu(P)=m!\vol{P}$ is an integer, where $\vol{P}$ is the usual
$m$-dimensional volume of~$P$.  A lattice simplex is called \emph{even} or \emph{odd} depending on
the parity of its normalized volume.  A triangulation~$K$ of a lattice polytope $P$ is \emph{dense}
if it uses all lattice points inside~$P$, that is, its vertex set is $P\cap\ZZ^m$.  In the case that
$K$ is additionally regular, say with lifting function~$\lambda$, we again write $P^\lambda$ for $K$
since it only depends on~$P$ and~$\lambda$.

Let $P^\lambda$ be an {rdf}-triangulation of~$P$, that is, $P^\lambda$ is regular, dense, and
foldable. In particular~$P^\lambda$ is a lattice triangulation.  Recall that $P^\lambda$ is foldable
if and only if its dual graph is bipartite.  Usually we refer to the two color classes as ``black''
and ``white''.  Then the \emph{signature} $\sigma(P^\lambda)$ of~$P^\lambda$ is defined as the
absolute value of the difference of the odd black and the odd white facets in~$P^\lambda$.  Note
that the even facets are not accounted for in any way.  Moreover, in the important special case
where $P^\lambda$ is \emph{unimodular}, that is, where all the facets have a normalized volume equal
to~$1$, all facets are odd.  For examples of unimodular triangulations of the $3$-cube with
signatures equal to~$0$ and~$2$ see Figure~\ref{fig:three_orderings}; note that all triangulations
of the $3$-cube without additional vertices are regular.

\begin{exmp}\label{exmp:undef}
  Dense and foldable triangulations do not exist
  for all lattice polytopes. For instance, in any dimension $m \ge 2$ there are lattice simplices of
  arbitrarily large volume which admit exactly one dense triangulation (which is regular), but which
  is not foldable.
  
  For $k \geq 1$ let $\Delta_2(k) = \conv\{(0,1),(1,0),(2k,2)\}$, a triangle with normalized
  volume $\nu(\Delta_2(k)) = 2k+1$. For $m \geq 3$ we define $\Delta_m(k)$ as the cone over
  $\Delta_{m-1}(k)$ with the $m$-th unit vector as its apex; this is an $m$-simplex with normalized
  volume $\nu(\Delta_m(k)) = \nu(\Delta_{m-1}(k)) = \ldots = 2k+1$.
  
  \begin{figure}[htbp]\centering
    \psfrag{(0,1)}[][]{$(1,0)$}
    \psfrag{(1,0)}[][]{$(0,1)$}
    \psfrag{(k,2)}[][]{$(2k,2)$}
    \includegraphics[width=.6\textwidth]{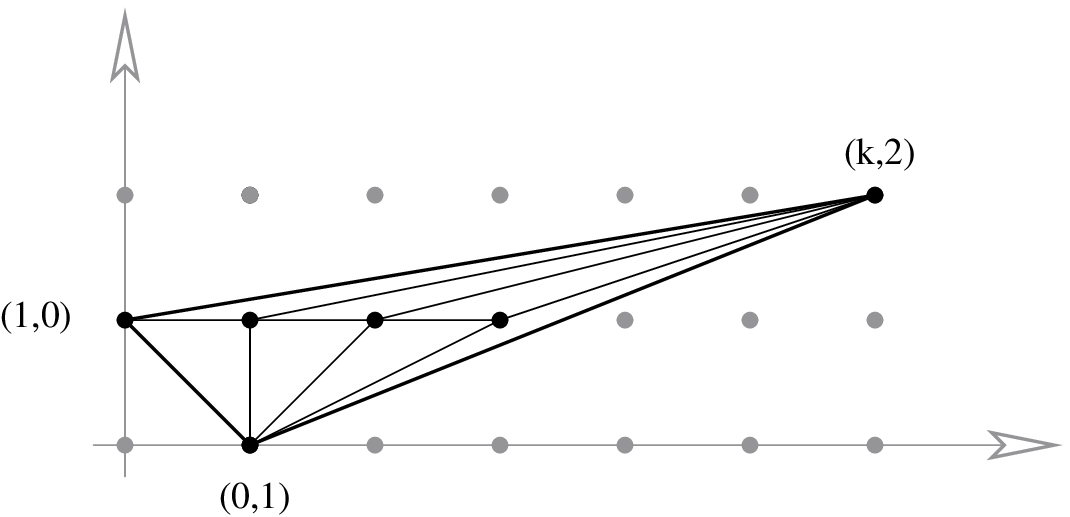}
  \end{figure}

  The interior point $(k,1) \in \Delta_2(k)$ is a degree-$3$-vertex in the unique (regular and)
  dense triangulation of~$\Delta_2(k)$, hence there is no dense and foldable triangulation
  of~$\Delta_2(k)$.  The cone over a triangulation~$K$ of~$\Delta_{m-1}(k)$ is foldable if and only
  if $K$ is foldable and any triangulation of~$\Delta_m(k)$ arises as a cone over a triangulation
  of~$\Delta_{m-1}(k)$.  Therefore there is no {rdf}-triangulation of~$\Delta_m(k)$ by
  induction.
\end{exmp}

\begin{exmp}[Signature of the Staircase Triangulation] \label{exmp:staircase_triang}
  Let $\Delta_m$ and $\Delta_n$ be odd simplices of dimension~$m$ and~$n$, respectively.  From the
  description of $\Gamma^*(\stc_{m,n})$ as the intersection of~$\mathcal{S}_{m,n}$
  with~$\mathcal{L}_m$ (see Proposition~\ref{prop:graph_as_intersection}) one can read off
  that~$\Gamma^*(\stc_{m,n})$ is bipartite and extract a recursive formulae for the signature
  of~$\stc_{m,n}$.  Remember that~$\stc_{m,n}$ is unimodular,
  hence~$\sigma_{m,0} = \sigma_{0,n} = 1$ and
\begin{align*}
  \sigma_{m,n} \; &= \; \abs{\; \sum_{i=0}^n (-1)^i \: \sigma_{m-1,i} \;}
  \; = \; \abs{\; \sum_{i=0}^{n-1}  (-1)^i \: \sigma_{m-1,i} \; + \; (-1)^n \: \sigma_{m-1,n} \;}\\
  &= \; \abs{\; \sigma_{m,n-1} \; + \; (-1)^n \: \sigma_{m-1,n} \;} \; = \; \sigma_{m,n-1} \; + \; (-1)^n \: \sigma_{m-1,n} \; .\\
\end{align*}

A careful inspection of the four cases arising from the two choices each for the parities of~$m$
and~$n$ gives the last equation.  This recursion then yields the explicit formulae for~$\sigma_{m,n}$
given by White~\cite{MR1840476} and stated in Proposition~\ref{prop:sign_m_n}.  Observe that
$\Delta_m\times\Delta_n$ is the order polytope of the poset of the disjoint union of a path of
length $m+1$ and a path of length $n+1$.  The staircase triangulation $\stc_{m,n}$ coincides with
the \emph{canonical triangulation} of the order polytope; see Soprunova and
Sottile~\cite[Section~4]{math.AG/0409504}.

\begin{prop} \label{prop:sign_m_n} The signature of the staircase triangulation of the product of
  two simplices of odd normalized volume is
 \[  
 \sigma_{2k,2l} = \binom{k+l}{k}\; , \quad
 \sigma_{2k,2l+1} = \binom{k+l}{k}\quad \text{and} \quad
 \sigma_{2k+1,2l+1} = 0 \; .
 \]
 If at least one of the simplices is even then this signature vanishes.
\end{prop}

\end{exmp}

Throughout the rest of the section let $P\subset\RR^m$ and $Q\subset\RR^n$ be an $m$- and
$n$-dimensional lattice polytopes, respectively.  Further we assume that there are
{rdf}-triangulations $P^\lambda$ and $Q^\mu$.  Suppose now that we have linear orderings~$O_P$ and~$O_Q$
of the vertex sets $V_P=P\cap\ZZ^m$ and $V_Q=Q\cap\ZZ^n$ such that the corresponding simplicial
product $P^\lambda\simpp Q^\mu$ is again foldable.  Note that such orderings always exist due to
Proposition~\ref{prop:foldable_product}.  By Proposition~\ref{prop:regular}, $P^\lambda\simpp Q^\mu$
is also regular and dense.

The rest of this section is devoted to computing the signature of $P^\lambda\simpp Q^\mu$.  The dual
graph~$\Gamma^*$ of the cell complex $P^\lambda\times Q^\mu$ is the product of the dual graphs of
$P^\lambda$ and $Q^\mu$.  Further the dual graph of the simplicial product $P^\lambda\simpp Q^\mu$
arises from $\Gamma^*$ by replacing each node by a copy of $\Gamma^*(\stc_{m,n})$ in a suitable way.

Recall that only odd simplices contribute to the signature.  Since the staircase triangulation is
unimodular for each facet $F$ of $\stc(f\times g)$ we have \mbox{$\nu(F)=\nu(f)\nu(g)$}.  Therefore we have
\begin{equation}\label{eq:signature}
  \sigma(P^\lambda\simpp Q^\mu)=\sigma_{m,n}\;\abs{\;\sum_{\text{$f\times g$ facet of~$P^\lambda\times Q^\mu$}\;}
    \delta(f,g)\;\overline{\nu}(f)\;\overline{\nu}(g)},
\end{equation}
where $\delta(f,g)=\pm1$ and $\overline{\nu}(h)=\nu(h)\mod 2$ denotes the parity of the normalized
volume of $h$.  So it remains to determine the sign $\delta(f,g)$.  This only depends on the
vertex orderings $O_P$ and $O_Q$.

As a \emph{point of reference} inside $\stc(f\times g)$ we choose the facet $F_0(f,g)$ corresponding to
the origin in the notation from Section~\ref{sec:prod_of_simp}; this corresponds to the
staircase $F_0=11\dots 100\dots 0$ which first goes all the way to the right and then all the way up
in Figure~\ref{fig:staircase}.  To determine the sign $\delta(f,g)$ amounts to determining the color
of the facet $F_0(f,g)$ in $P^\lambda\simpp Q^\mu$. 

We first consider the case where $P^\lambda$ is a lattice $m$-simplex~$\Delta_m$ (without interior
lattice points) and~$Q^\mu$
consists of two neighboring $n$-simplices (without interior lattice points), that is,~$Q^\mu$ is the
{rdf}-triangulation~$B_n$ of the bipyramid over the $(n-1)$-simplex from Example~\ref{exmp:bipyramid}.
Note that~$\Delta_m$ is an {rdf}-triangulation of itself.  Further,
the signature of~$\Delta_m$ vanishes if the normalized volume of~$\Delta_m$ is even and equals~$1$
otherwise.

\begin{lem}\label{lem:D_X_Bn}
  The simplicial product $\Delta_m \simpp B_n$ is an {rdf}-triangulation of the product
  of~$\Delta_m$ and a lattice bipyramid over the $(n-1)$-simplex with signature
  \[\sigma(\Delta_m \simpp B_n)\;=\;
  \begin{cases}
    \sigma_{m,n}\;\sigma(\Delta_m)\;\sigma(B_n) & \text{\begin{tabular}{l}
                                                           if the vertex ordering on $B_n$ is\\
                                                           color consecutive or if $m$ is even,
                                                         \end{tabular}
                                                       }\\
    \sigma_{m,n}\;\sigma(\Delta_m)\;\omega  & \text{\begin{tabular}{l}
                                                      if the vertex ordering on $B_n$\\
                                                      is symmetric and $m$ is odd.
                                                    \end{tabular}
                                                  }
    \end{cases}
  \]
  Here $\omega\in\{0,1,2\}$ counts the number of odd simplices in $B_n$.
\end{lem}

\begin{proof}
  It is a consequence of Propositions~\ref{prop:foldable_product} and~\ref{prop:regular} that
  $\Delta_m \simpp B_n$ is an {rdf}-triangulation.
  
  Let $g$ and $g'$ be the two facets of $B_n$.  In both cases we get a contribution of
  $\delta(\Delta_m,g)\;\sigma_{m,n}\;\sigma(\Delta_m)$ to $\sigma(\Delta_m \simpp B_n)$ if~$g$ is
  odd, and similarly a contribution of $\delta(\Delta_m,g')\;\sigma_{m,n}\;\sigma(\Delta_m)$ to
  $\sigma(\Delta_m \simpp B_n)$ if $g'$ is odd; see Equation~\eqref{eq:signature}.
  
  It remains to compare $\delta(\Delta_m,g)$ and $\delta(\Delta_m,g')$, which depends on the
  vertex ordering of $B_n$. We have $\delta(\Delta_m,g)=-\delta(\Delta_m,g')$ if and only if
  $F_0(\Delta_m,g)$ and $F_0(\Delta_m,g')$ are colored differently which in turn holds if and only
  if the distance between $F_0(\Delta_m,g)$ and $F_0(\Delta_m,g')$ in $\Gamma^*(\Delta_m \simpp
  B_n)$ is odd.

  Since $\Gamma^*(\Delta_m \simpp B_n)$ is bipartite, each path from $F_0(\Delta_m,g)$ to
  $F_0(\Delta_m,g')$ has the same parity, and we may choose any path to determine the parity of the
  distance. Let $\tilde{F}_0(\Delta_m,g) \in \stc(\Delta_m\times g)$ and $\tilde{F}_0(\Delta_m,g')
  \in \stc(\Delta_m\times g')$ be neighboring facets. Then the distance between $F_0(\Delta_m,g)$
  and $F_0(\Delta_m,g')$ is odd if and only if the distance between $F_0(\Delta_m,g)$ and
  $\tilde{F}_0(\Delta_m,g)$ has the same parity as the distance between $F_0(\Delta_m,g')$ and
  $\tilde{F}_0(\Delta_m,g')$ (keep in mind that the distance between $\tilde{F}_0(\Delta_m,g)$ and
  $\tilde{F}_0(\Delta_m,g')$ is~$1$).
  
  We first consider the case where the vertex ordering of $B_n$ is color consecutive.  Let~$c$ be
  the color of the unique vertex $v\in g\setminus g'$ (which is the same as the color of the unique
  vertex $v'\in g'\setminus g$).  All columns in the lattice grid defining $\Delta_m \simpp B_n$
  corresponding to vertices colored~$c$ are consecutive and hence $v$ and~$v'$ follow one after
  another in the vertex ordering of~$B_n$.  We distinguish the two cases where~$v$ and~$v'$ appear
  somewhere in the middle or at the beginning of the vertex ordering of~$B_n$ and where~$v$ and~$v'$
  appear at the end of the vertex ordering; see Figure~\ref{fig:consecutive_bipyramid}.  In the
  first case we may choose $F_0(\Delta_m,g)=\tilde{F}_0(\Delta,g)$ and
  $F_0(\Delta_m,g')=\tilde{F}_0(\Delta_m,g')$ and the distance between $F_0(\Delta_m,g)$ and
  $F_0(\Delta,g')$ is~$1$.  In the second case the distance between $F_0(\Delta_m,g)$ and
  $\tilde{F}_0(\Delta_m,g)$ equals the distance between $F_0(\Delta_m,g')$ and
  $\tilde{F}_0(\Delta_m,g')$.  Therefore we obtain $\delta(\Delta_m,g)=-\delta(\Delta_m,g')$ in
  the color consecutive case.

  \begin{figure}[htbp]\centering
    \begin{overpic}[width=.43\textwidth]{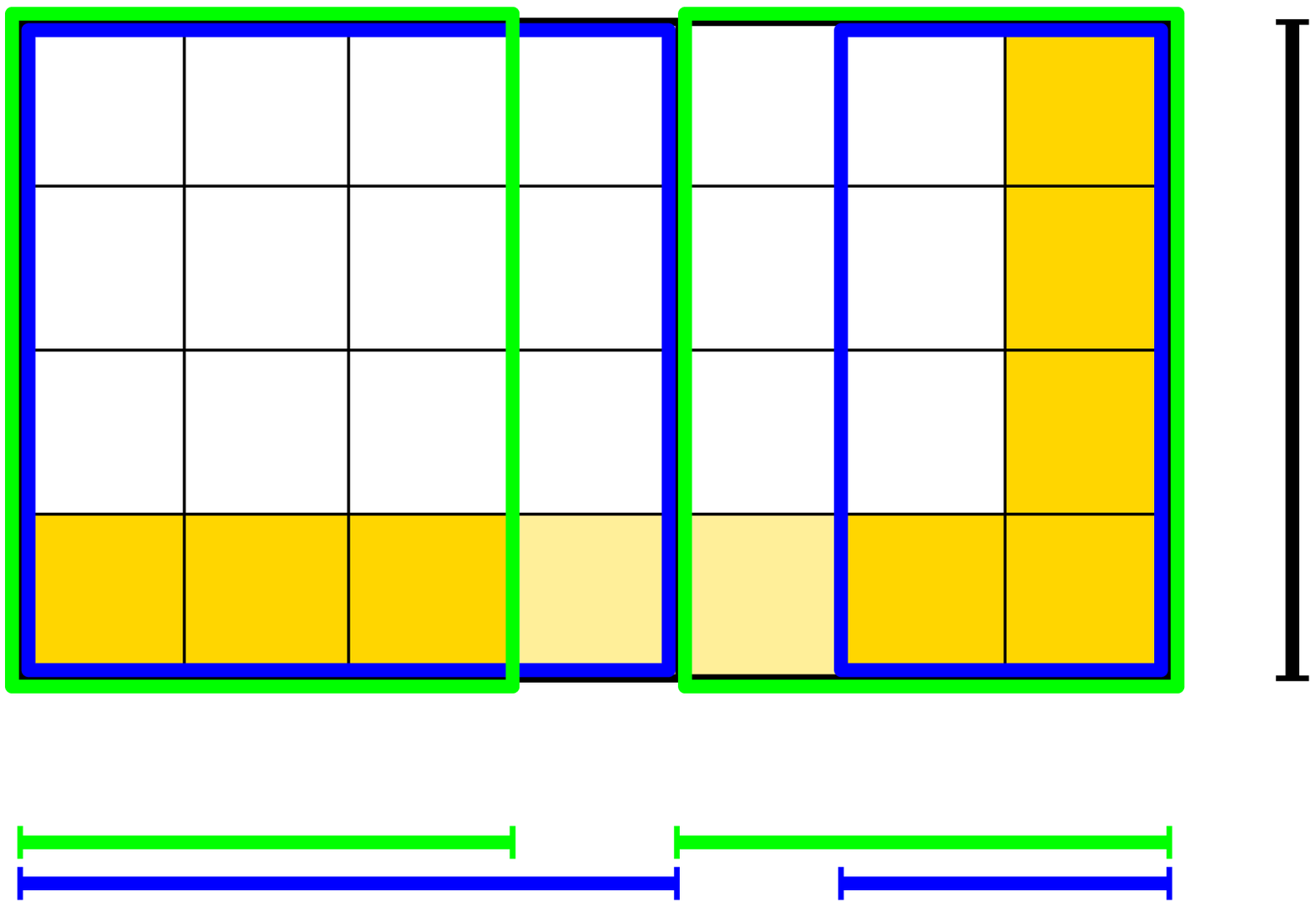}
      \put(97,50){$\Delta_m$}
      \put(14,4){$g \cap g'$}
      \put(67,4){$g \cap g'$}
      \put(42,4){$v$}
      \put(54,4){$v'$}
    \end{overpic}
    \hspace{0.4cm}
    \begin{overpic}[width=.43\textwidth]{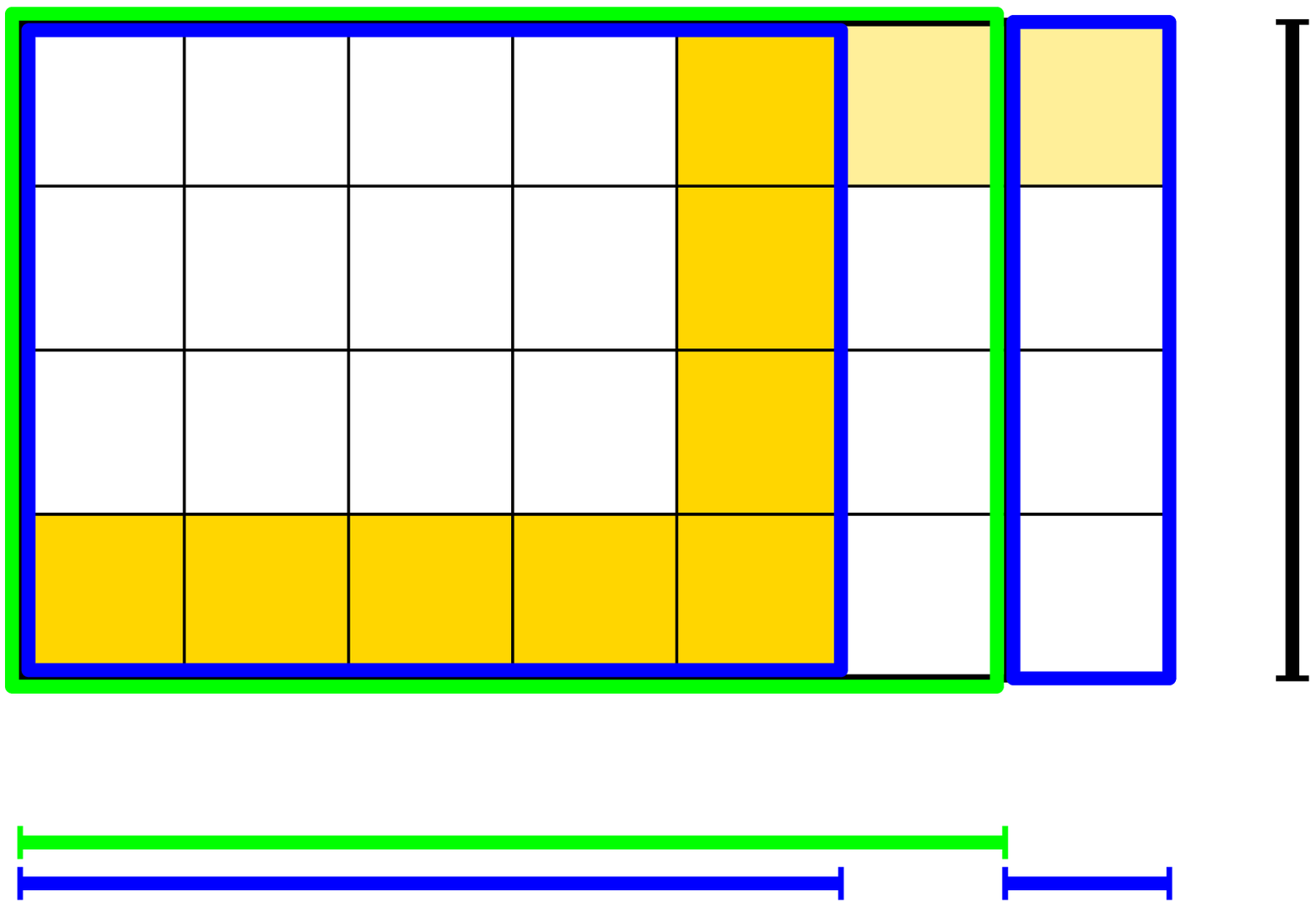}
      \put(97,50){$\Delta_m$}
      \put(26,4){$g \cap g'$}
      \put(65,4){$v$}
      \put(78,4){$v'$}
    \end{overpic}
    \caption{Distance of the facets of reference $F_0(\Delta_m,g)$ and $F_0(\Delta_m,g')$ in
      $\Gamma^*(\Delta_m \simpp B_n)$ for color consecutive orderings of $B_n$.  The facets
      $\tilde{F}_0(\Delta_m,g)$ and $\tilde{F}_0(\Delta_m,g')$ and their intersection is shaded.  On
      the left the two apices $v,v'$ occur somewhere in the middle or at the beginning of the vertex
      ordering of~$B_n$, on the right at the end.\label{fig:consecutive_bipyramid}}

    \vspace{.4cm}
    \begin{overpic}[width=.43\textwidth]{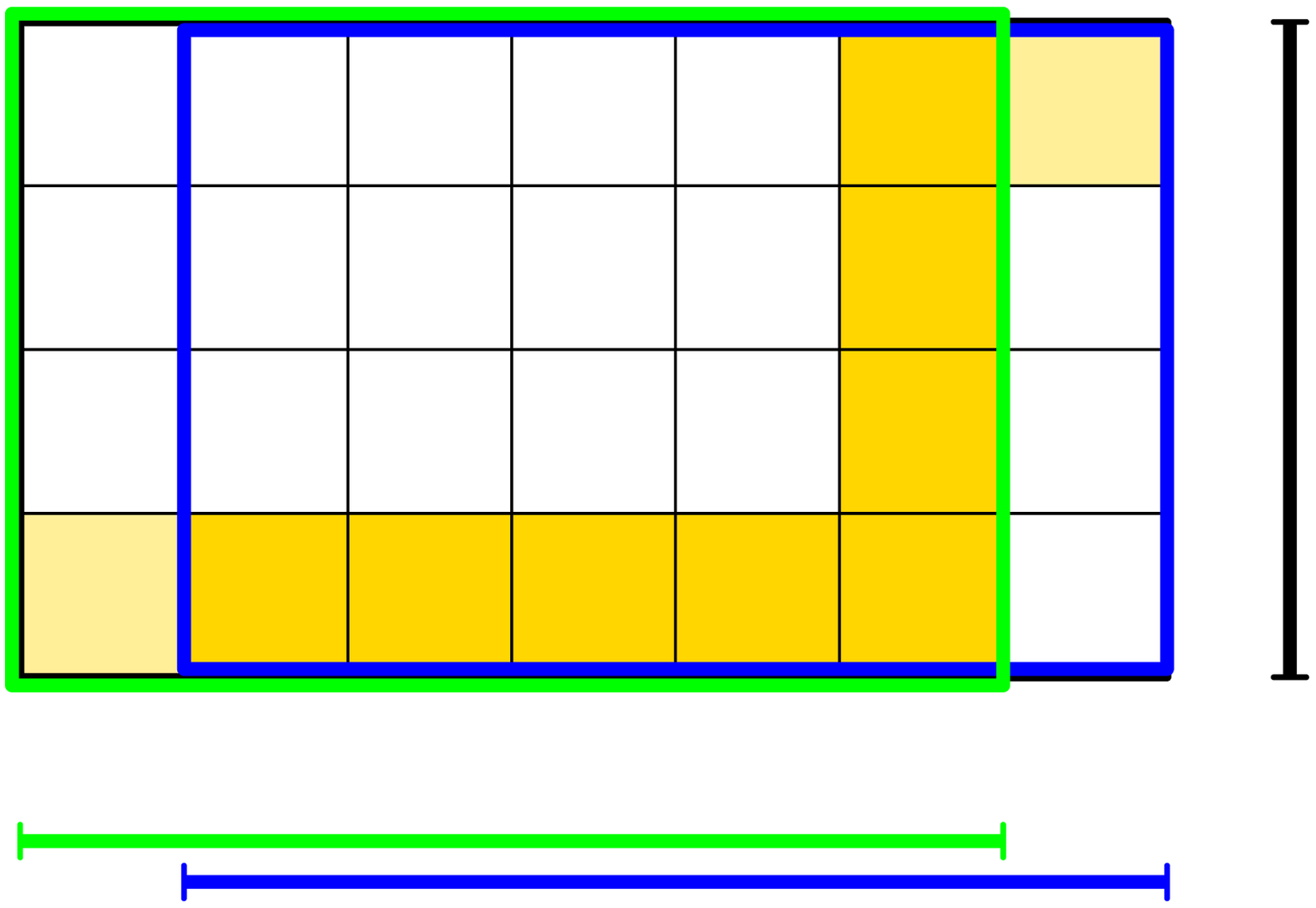}
      \put(97,50){$\Delta_m$}
      \put(38,4){$g \cap g'$}
      \put(6,4){$v$}
      \put(78,4){$v'$}
    \end{overpic}
    \caption{Distance of the facets of reference $F_0(\Delta_m,g)$ and $F_0(\Delta_m,g')$ in
      $\Gamma^*(\Delta_m \simpp B_n)$ for the symmetric ordering of the vertices of $B_n$.  The
      facets $\tilde{F}_0(\Delta_m,g)$ and $\tilde{F}_0(\Delta_m,g')$ and their intersection is
      shaded.\label{fig:non_consecutive_ordering}}
  \end{figure}

  Let the vertex ordering on $B_n$ be symmetric.  We have $F_0(\Delta_m,g) =
  \tilde{F}_0(\Delta_m,g)$ and the distance of $F_0(\Delta_m,g')$ and $\tilde{F}_0(\Delta_m,g')$ is
  $m$, hence $\delta(\Delta_m,g)=-\delta(\Delta_m,g')$ if and only if $m$ is even; see
  Figure~\ref{fig:non_consecutive_ordering}.

  We refer to Figure~\ref{fig:three_orderings} for an
  example of three triangulations of $[0,1]\times B_2$ resulting from different vertex orders
  of~$B_2$.
\end{proof}

\begin{thm}[Combinatorial Product Theorem] \label{thm:comb_product_thm}
  Let $P^\lambda$ and $Q^\mu$ be {rdf}-triangulations of an $m$-dimensional lattice
  polytope $P\subset\RR^m$ and an $n$-di\-men\-si\-o\-nal lattice polytope $Q\subset\RR^n$, respectively.
  For color consecutive vertex orderings $O_P$ and $O_Q$ the simplicial product $P^\lambda\simpp
  Q^\mu$ is an {rdf}-triangulation of the polytope $P\times Q$ with signature
  \[\sigma(P^\lambda\simpp Q^\mu)=\sigma_{m,n}\;\sigma(P^\lambda)\;\sigma(Q^\mu)\;.\]
\end{thm}

\begin{proof}
  Again, by Propositions~\ref{prop:foldable_product} and~\ref{prop:regular}, $P^\lambda\simpp Q^\mu$
  is an {rdf}-triangulation.
  
  \enlargethispage{\baselineskip} Let~$f,f' \in P^\lambda$ and $g,g' \in Q^\mu$ be facets such that
  $f \times g$ and $f' \times g'$ are neighboring cells of $P^\lambda \times Q^\mu$. We may assume
  that $f=f'$ and $g \cap g'$ is a ridge.  Hence $g\cup g'$ is a bipyramid over the common ridge
  $g\cap g'$.  Applying Lemma~\ref{lem:D_X_Bn} to $f\simpp(g\cup g')$ yields $\delta(f,g)=
  -\delta(f,g')$, and we may label the cells of $P^\lambda \times Q^\mu$ with $\delta(f,g)$ by
  assigning~$+1$ (black) and~$-1$ (white) according to the bipartition of the dual graph
  $\Gamma^*(P^\lambda \times Q^\mu)$ of $P^\lambda \times Q^\mu$.
  
  \enlargethispage{2\baselineskip}
  We may think of $\Gamma^*(P^\lambda \times Q^\mu)$ as a copy of $\Gamma^*(P^\lambda)$ for each node of $\Gamma^*(Q^\mu)$.
  Each copy of $\Gamma^*(P^\lambda)$ may be $2$-colored using the bipartition of
  $\Gamma^*(P^\lambda)$, but we must use the inverse coloring for a copy of~$\Gamma^*(P^\lambda)$ if
  the corresponding node of~$\Gamma^*(Q^\mu)$ is colored white. Therefore a node $f \times g$ of
  $\Gamma^*(P^\lambda \times Q^\mu)$ is labeled~$+1$ if and only if the facets $f \in P^\lambda$ and $g \in Q^\mu$
  are colored the same, and using Equation~\eqref{eq:signature} we have
    \begin{align*}
    \sigma(P^\lambda \simpp Q^\mu)
    =\;&\sigma_{m,n} \left|\, \sum_{\text{$f\in P^\lambda$
          black}}\hspace{-.15cm}\Big(\,\overline{\nu}(f)\hspace{-.15cm}\sum_{\text{$g\in Q^\mu$
          black}}\hspace{-.2cm}\overline{\nu}(g)\,\Big) \right.
    \;+ \hspace{-.3cm}\sum_{\text{$f\in P^\lambda$ white}}\hspace{-.15cm}\Big(\,\overline{\nu}(f)\hspace{-.15cm}\sum_{\text{$g\in Q^\mu$
          white}}\hspace{-.2cm}\overline{\nu}(g)\,\Big)\\
    &- \hspace{-.3cm}\sum_{\text{$f\in P^\lambda$
          black}}\hspace{-.15cm}\Big(\,\overline{\nu}(f)\hspace{-.15cm}\sum_{\text{$g\in Q^\mu$
          white}}\hspace{-.2cm}\overline{\nu}(g)\,\Big)
    \;- \hspace{-.3cm}\left. \sum_{\text{$f\in P^\lambda$ white}}\hspace{-.15cm}\Big(\,\overline{\nu}(f)\hspace{-.15cm}\sum_{\text{$g\in Q^\mu$
          black}}\hspace{-.2cm}\overline{\nu}(g)\,\Big) \,\right|\\
    =\; &\sigma_{m,n} \left| \sum_{\text{$f\in P^\lambda$ black}}\hspace{-.2cm}\overline{\nu}(f) \;- \hspace{-.3cm}
        \sum_{\text{$f\in P^\lambda$ white}}\hspace{-.2cm}\overline{\nu}(f)\right| \,
      \left|\sum_{\text{$g\in Q^\mu$ black}}\hspace{-.2cm}\overline{\nu}(g) \;- \hspace{-.3cm} \sum_{\text{$g\in Q^\mu$
            white}}\hspace{-.2cm}\overline{\nu}(g)\right|\\
    =\; &\sigma_{m,n} \; \sigma(P^\lambda) \; \sigma(Q^\mu)\;.
    \end{align*}
\end{proof}

Finally we consider the case where~$Q^\mu$ is the {rdf}-triangulation~$B_n$ of the
bipyramid over the $(n-1)$-simplex from Example~\ref{exmp:bipyramid}.  While this seems to cover a
very special case only, the result is instrumental for the construction of triangulations of the
$d$-cube with non-trivial signature in Section~\ref{sec:cubes}.

\begin{prop} \label{prop:sign_of_bipyramid}
  Let $P^\lambda$ be an {rdf}-triangulation of an $m$-di\-men\-si\-o\-nal lattice polytope
  $P\subset\RR^m$ with a color consecutive ordering on its vertex set
  $V_P=P\cap\ZZ^m$.
  Then $P^\lambda\simpp B_n$ is an {rdf}-triangulation of the product of $P$ with a lattice bipyramid
  over the $(n-1)$-simplex with signature
  \[\sigma(P^\lambda\simpp B_n)\;=\;
  \begin{cases}
    \sigma_{m,n}\;\sigma(P^\lambda)\;\sigma(B_n) & \text{\begin{tabular}{l}
                                                           if the vertex ordering on $B_n$ is\\
                                                           color consecutive or if $m$ is even,
                                                         \end{tabular}
                                                       }\\
    \sigma_{m,n}\;\sigma(P^\lambda)\;\omega & \text{\begin{tabular}{l}
                                                      if the vertex ordering on $B_n$\\
                                                      is symmetric and $m$ is odd.
                                                    \end{tabular}
                                                  }
  \end{cases}
  \]
  Here $\omega\in\{0,1,2\}$ counts the number of odd simplices in $B_n$.
\end{prop}

One can show that for other vertex orderings of $B_n$ the simplicial product $P^\lambda\simpp B_n$ is not
foldable.  In this sense the two cases listed exhaust all the possibilities.

\begin{proof}
  Propositions~\ref{prop:foldable_product} and~\ref{prop:regular} ensure that
  $P^\lambda\simpp Q^\mu$ is an {rdf}-triangulation.
  Let $g$ and $g'$ be the two facets of $B_n$, and let us think of $P^\lambda \times B_n$ as the
  union of two copies of $P^\lambda \times \Delta_n$, which we denote as $P^\lambda \times g$ and
  $P^\lambda \times g'$.  Further let $f \in P^\lambda$ be an arbitrary but fixed facet. We get a
  contribution of $\delta(f,g)\;\sigma(P^\lambda)\;\sigma_{m,n}$ to $\sigma(P^\lambda \simpp B_n)$
  if~$g$ is odd by Theorem~\ref{thm:comb_product_thm}. Similarly we get a contribution of
  $\delta(f,g')\;\sigma(P^\lambda)\;\sigma_{m,n}$ to $\sigma(P^\lambda \simpp B_n)$ if~$g'$ is odd.
  It remains to compare $\delta(f,g)$ and $\delta(f,g')$.  The simplicial product $f\simpp (g\cup
  g')$ is a triangulation of the product of an $m$-simplex and~$B_n$ and by Lemma~\ref{lem:D_X_Bn}
  we have $\delta(f,g)=-\delta(f,g')$ in the first and $\delta(f,g)=\delta(f,g')$ in the second case.
\end{proof}

A referee suggested the following generalization of Proposition~\ref{prop:sign_of_bipyramid}, which
we state without a proof.  Let~$P^\lambda$ and~$Q^\mu$ be {rdf}-triangulations of the full
dimensional lattice polytopes $P\subset\RR^m$ and $Q\subset\RR^n$, respectively. Further let the
vertices of~$P^\lambda$ be ordered color consecutive, and let the vertices of~$Q^\mu$ be partitioned
into subsets $V_0,V_1,\dots,V_n$ according to their colors. An \emph{almost color consecutive}
ordering of the vertices of~$Q^\mu$ is obtained by splitting~$V_0$ into two subsets~$V_0'$
and~$V_0''$ and taking any vertex ordering compatible with $V_0'<V_1<\dots<V_n<V_0''$. The vertex
sets~$V_0'$ and~$V_0''$ induce a bipartition on the facets of $Q^\mu$ denoted by~$L'$ and~$L''$, and
let the facets of~$L'$, respectively~$L''$ be colored ``black'' and ``white'' according to the
coloring of the facets of~$Q^\mu$ (neither~$L'$ nor~$L''$ is strongly connected in general). Finally
we set the \emph{signed signature}~$\tilde{\sigma}(L)$ of a geometric simplicial complex~$L$ with
facets colored ``black'' and ``white'' as the number of odd ``black'' facets minus the number of odd
``white'' facets.

\begin{prop}\label{prop:sign_of_bipyramid2}
  The simplicial product $P^\lambda\simpp Q^\mu$ (with respect to the color consecutive vertex
  ordering of~$P^\lambda$ and the almost color consecutive vertex ordering of~$Q^\mu$) is a
  {rdf}-triangulation of $P\times Q$ with signature
  \[ 
  \sigma(P^\lambda\simpp Q^\mu) = \begin{cases} 
    \sigma_{m,n}\;\sigma(P^\lambda)\;\sigma(Q^\mu)
    & \text{if $m$ is even,}\\
    \sigma_{m,n}\;\sigma(P^\lambda)\;|\tilde{\sigma}(L')-\tilde{\sigma}(L'')|
    & \text{if $m$ is odd.}
  \end{cases}
  \]
\end{prop}

\section{Lower Bounds for the Number of Real Roots of Polynomial Systems}

Triangulations which are regular, dense, and foldable are interesting since they yield
non-trivial lower bounds for the number of real roots of associated polynomial systems, provided
that a number of additional geometric conditions are met.  To discuss these issues we first review
the construction of Soprunova and Sottile~\cite{math.AG/0409504}.

\subsection{Triangulations and Lower Bounds}\label{sec:real_roots}

Let~$P\subset \RR_{\ge 0}^m$ be a lattice $m$-polytope contained in the positive orthant, and let
$\lambda:P \cap\ZZ^m\to\RR$ be a lifting function such that the induced triangulation~$P^\lambda$
is an {rdf}-triangulation.  Further let the vertices $P \cap\ZZ^m$ of~$P^\lambda$ be colored by the map $c
: P \cap \ZZ^m \to [m+1]$.  We define the \emph{coefficient polynomial} $F_{P^\lambda,i,s}\in
\RR[t_1,\dots,t_m]$ of a color $i$ and an additional parameter $s\in (0,1]$ as
\begin{equation} \label{eq:coeff_poly}
F_{P^\lambda,i,s}(t)  \;=\; \sum_{v \: \in \: c^{-1}(i)} s^{\lambda(v)} \; t^v,
\end{equation}
where $t = (t_1, \dots, t_m)$ and $t^v = t_1^{v_1}\dots t_m^{v_m}$.  Choosing a real number
$a_i$ for each color $i\in[m+1]$ defines a \emph{Wronski polynomial}
\[
\mathcal{F}_{P^\lambda,s}(t) \;=\;  a_0 F_{P^\lambda,0,s}(t) \; + \; a_1 F_{P^\lambda,1,s}(t)
\; + \; \ldots \; + \; a_m F_{P^\lambda,m,s}(t) \; \in \RR[t_1,\dots,t_m]\;,
\]
for fixed $s \in (0,1]$.  A \emph{Wronski system} associated with $P^\lambda$ is a sparse system of
$m$ Wronski polynomials which is \emph{generic} in the sense that it attains Kushnirenko's
bound~\cite{kushnirenko:newton_polyhedron}, that is, it has exactly $\nu(P)$ distinct complex
solutions.

Let $M=\abs{P \cap \ZZ^m}$ denote the number of integer points in $P$ and let $\CC\PP^{M-1}$ be the
complex projective space with coordinates $\smallSetOf{x_v}{v \in P \cap \ZZ^m}$.  The toric
projective variety $X_P \subset \CC\PP^{M-1}$ parameterized by the monomials $\smallSetOf{t^v}{v \in
  P \cap \ZZ^m}$ is given by the closure of the image of the map
\begin{equation}\label{eq:phi}
\phi_P : (\CC^{\times})^m \to \CC\PP^{M-1} : t \mapsto [ t^v \; | \; v \in P \cap \ZZ^m ]\;,
\end{equation}
where $[t^{v_1},\dots,t^{v_m}]$ is a point in $\CC\PP^{M-1}$ written in homogeneous coordinates.
Via~$\phi_P$ a Wronski system on $(\CC^{\times})^m$ corresponds to a system of $m$ linear equations on
the toric variety $X_P\subset \CC\PP^{M-1}$.

Let $Y_P = X_P \cap \RR\PP^{M-1}$ be the real points of the variety $X_P$.  For $s \in (0,1]$ the
\emph{$s$-deformation}~$s.Y_P$ is obtained as the closure of the image of the deformed map
\[
s.\phi_P : (\CC^{\times})^m \to \CC\PP^{M-1} : t \mapsto [ s^{\lambda(v)}\;t^v \; | \; v \in P \cap \ZZ^m ]
\]
intersected with~$\RR\PP^{M-1}$. The $s$-deformation~$s.Y_P$ interpolates between $Y_P=1.Y_P$ and
its homotopic image $0.Y_P$, which is defined as the initial variety $\Initial_\lambda(Y_P)$; the
whole family $\smallSetOf{s.Y_P}{s\in[0,1]}$ is called the \emph{toric degeneration} of~$Y_P$; for
the details see~\cite[Section~3]{math.AG/0409504}.  A Wronski polynomial corresponds to the image of
$s.Y_P$ under the linear \emph{Wronski projection}
\[
\pi_E :\;
\begin{aligned}
  \CC\PP^{M-1} \setminus E         \,&\to\, \CC\PP^m\\
  [x_v \; | \; v \in P \cap \ZZ^m] \,&\mapsto\, [ \sum_{v \: \in \: c^{-1}(i)} x_v \; | \; i =
  0,1,\dots,m \; ]
\end{aligned}
\]
with \emph{center}
\[
E = \SetOf{x \in \CC\PP^{M-1}}{\sum_{v \: \in \: c^{-1}(i)} x_v = 0 \quad \text{for $i = 0,1,\dots,m$}}.
\]

The toric degeneration \emph{meets} the center of the projection $\pi_E$ if there are $s
\in (0,1]$ and $t \in \RR^m$ such that
\[
F_{P^\lambda,0,s}(t) \; = \; F_{P^\lambda,1,s}(t) \; = \; \dots \; = \; F_{P^\lambda,m,s}(t) \; = \; 0 \; .
\]

The sphere $\Sph^{M-1}$ is a double cover of $\RR\PP^{M-1}$. Let $Y_P^+ \subset \Sph^{M-1}$ be the
pre-image of~$Y_P$ under the covering map.  Note that $Y_P^+$ is not necessarily smooth nor
connected.  Nonetheless, its orientability is well defined.  The following theorem is a slightly
simplified version of what is proved in~\cite{math.AG/0409504}.

\begin{thm}[Soprunova \& Sottile] \label{thm:Soprunova_Sottile}
  Let $P \subset \RR_{\ge0}^m$ be a non-negative lattice $m$-polytope such that $Y_P^+$ is oriented,
  and let $P^\lambda$ be an {rdf}-triangulation of~$P$ induced by the lifting function~$\lambda$.
  Suppose that there is a number $s_0 \in (0,1]$ such that the $s$-deformation $s.Y_P$ does not meet
  the center of the Wronski projection~$\pi_E$ for all $s \in (0,s_0]$ and all $t\in\RR^m$.  Then
  for all $s \in (0,s_0]$ the number of real solutions of any associated Wronski system in
  $\RR[t_1,\dots,t_m]$ is bounded from below by the signature $\sigma(P^\lambda)$.
\end{thm}

In general, it seems difficult to decide the orientability of $Y_P^+$.  To this end Soprunova and
Sottile suggest to consider the following sufficient condition: Let $(A,b)$ be an integral facet
description of $P=\smallSetOf{x\in\RR^m}{Ax+b\geq 0}$ such that the $i$-th row of the matrix~$A$ is
the unique inward pointing primitive normal vector of the $i$-th facet of~$P$.  This way, up to a
re-ordering of the facets,~$A$ and~$b$ are uniquely determined.  Denote by~$\Lambda_A$ the lattice
spanned by the columns of~$A$.  Suppose that the lattice spanned by $P \cap \ZZ^m$ has odd index in
$\ZZ^m$ and that~$\Lambda_A$ has odd index in its saturation $\Lambda_A \otimes_{\ZZ} \QQ$, that
is,~$A$ has a maximal minor~$\tilde{A}$ with $\det \tilde{A}$ odd.  If these two parity conditions
are satisfied and if, additionally, there is a vector~$v$ with only odd entries in the integer
column span of $(A,b)$ then Soprunova and Sottile call the double cover $Y_P^+$ \emph{Cox-oriented}.

We call the {rdf}-triangulation $P^\lambda$ \emph{nice} for the value $s_0$ if all the
conditions of Theorem~\ref{thm:Soprunova_Sottile} are satisfied.  Note that the (Cox-)orientability
of~$Y_P^+$ solely depends on the polytope~$P$.

\begin{exmp}
  The unique {rdf}-triangulation of the line segment $[k,l]$, where $0\le k<l$, is nice for
  $s_0=1$ (and any lifting function) if and only if $k=0$.  We have $\sigma([0,l])\in\{0,1\}$
  depending on $l$ being even or odd.  This is a sharp lower bound for the number of real roots in
  the one-dimensional case.
\end{exmp}

\begin{exmp}
  The staircase triangulation of $\Delta_m \times \Delta_n$ is nice for $s_0=1$.  This is true at
  least if one of the two vertices whose color occurs only once is located at the origin.
\end{exmp}

\begin{exmp}
  Let $P^\lambda$ be an {rdf}-triangulation of a lattice polytope $P\subset\RR_{\ge
    0}^m$, and let $Y_P^+$ be Cox-oriented.  The cone $0*P^\lambda$ of the triangulation~$P^\lambda$
  (embedded into $\RR^{m+1}$ via the map $(v_1,\dots,v_m)\mapsto(1,v_1,\dots,v_m)$) with apex
  $0\in\RR^{m+1}$ is nice for $s_0=1$.  The signature of $0*P^\lambda$ equals the signature of
  $P^\lambda$.
\end{exmp}

\subsection{Products of Toric Varieties}
\label{sec:poly_sys_from_prod}

Let us consider the \emph{Segre embedding}
\[
  \iota  : \,
  \begin{aligned}
    \CC\PP^{M-1}\times\CC\PP^{N-1}\,&\to\,\CC\PP^{MN-1}\\
    ([x_1,\dots,x_M],[y_1,\dots,y_N])\,&\mapsto\,[x_1y_1,\dots,x_iy_j,\dots,x_My_N]\; ,
  \end{aligned}
\]
which is the tensor product.  The restriction $\iota:\RR\PP^{M-1}\times\RR\PP^{N-1}\to\RR\PP^{MN-1}$
lifts to the double covers $\iota:\Sph^{M-1}\times\Sph^{N-1}\to\Sph^{MN-1}$.

\begin{prop}\label{prop:product_of_toric_varieties}
  Let $P$ be an $m$-dimensional lattice polytope with~$M$ lattice points, and let~$Q$ be an
  $n$-dimensional lattice polytope with~$N$ lattice points. Then we have
  \[
  \iota(Y_P \times Y_Q)\,=\,Y_{P\times Q} \quad\text{and}\quad \iota(Y_P^+ \times Y_Q^+)\,=\,Y_{P^+\times Q^+} \; .
  \]
\end{prop}

\begin{proof}
  Let $\phi_P:(\CC^{\times})^m \to \CC\PP^{M-1}$ denote the map in Equation~\eqref{eq:phi} which
  defines the toric variety~$X_P$.  Observe that $\phi_{P\times Q}=\iota\circ(\phi_P,\phi_Q)$.  This
  readily implies $\iota(X_P \times X_Q)=X_{P\times Q}$ and also $\iota(Y_P \times Y_Q)=Y_{P\times
    Q}$.  Now $\iota(Y_P^+ \times Y_Q^+)=Y_{P^+\times Q^+}$ follows since the map $\iota$ lifts to
  the coverings.
\end{proof}

\begin{cor} \label{cor:oriented}
  Let $P$ and $Q$ be lattice polytopes such that $Y_P^+$ and $Y_Q^+$ are oriented. Then $Y_{P \times Q}^+$
  is oriented.
\end{cor}

\begin{proof}
  The orientability of $Y_{P\times Q}^+$ depends on the orientability of its smooth part, which is
  the $\iota$-image of the product of the smooth parts of $Y_P^+$ and $Y_Q^+$.  The product of
  orientable manifolds is orientable.
\end{proof}

\begin{rem}
  As a further consequence, if $Y_P^+$ and $Y_Q^+$ are Cox-oriented, then $Y_{P \times Q}^+$ is
  oriented.  However, $Y_{P \times Q}^+$ does not have to be Cox-oriented itself.  For an example
  consider products $\Delta_m\times\Delta_n$ of standard simplices for $m$ even and $n$ odd.
\end{rem}

The question under which conditions the toric degeneration of $Y_{P\times Q}$ meets the center of the
Wronski projection is a little harder to answer.  The lifting function $\omega$ determines the
triangulation of $P \times Q$ and we write $(P \times Q)^\omega = P^\lambda \simpp Q^\mu$ if we want
to emphasize the particular lifting function~$\omega$ defined in Equation~\eqref{eq:omega1}.  Recall
that a vertex $(v,w)$ of $(P \times Q)^\omega$ is colored $k = c_{P^\lambda}(v) + c_{Q^\mu}(w)$
where $c_{P^\lambda}:P\cap\ZZ^m\to[m+1]$ and $c_{Q^\mu}:Q\cap\ZZ^n\to[n+1]$ denote the coloring
maps; see Proposition~\ref{prop:foldable_product}.  Therefore for $s \in (0,1]$ the coefficient
polynomial (Equation~\eqref{eq:coeff_poly}) of $(P \times Q)^\omega$ for $k\in[m+n+1]$ has the form
\begin{align*}
  F_{(P \times Q)^\omega,k,s}(t) \; &= \; \sum_{c_{P^\lambda}(v) + c_{Q^\mu}(w) = k}
  s^{\lambda(v)+\mu(w)+\epsilon(v,w)} \; t^{(v,w)}\\
  &= \; \sum_{c_{P^\lambda}(v) + c_{Q^\mu}(w) = k}
  s^{\lambda(v)}(t_1,\dots,t_m)^v \; s^{\mu(w)}(t_{m+1},\dots,t_{m+n})^w \; s^{\epsilon(v,w)} \; .
\end{align*}
As in Example~\ref{exmp:lifting_function_stc} we may choose the same perturbation $\epsilon(i,j) =
\epsilon \: 2^{(n+1)i+(n-j)}$ (for sufficiently small $\epsilon > 0$) for all vertices $(v,w)$ with
$c_{P^\lambda}(v)=i$ and $c_{Q^\mu}(w)=j$ if we choose color consecutive orderings of the vertices
of $P^\lambda$ and $Q^\mu$; see Equation~\eqref{eq:omega2}.  Summing over all colors
$i$ of $P^\lambda$ and all colors $j$ of $Q^\mu$ with $i+j=k$ yields
\begin{equation} \label{eq:prod_coeff_poly}
  F_{(P \times Q)^\omega,k,s} \; = \; \sum_{i+j=k} F_{P^\lambda,i,s} \; F_{Q^\mu,j,s} \; s^{\epsilon(i,j)} \;.
\end{equation}

The $s$-degeneration $s.Y_P$ meets the center of the Wronski projection in the points
\[
V_s({P^\lambda}) = \SetOf{t \in \RR^{m}}{F_{P^\lambda,i,s}(t)=0 \; \text{for all}
  \; i\in[m+1]} \; ,
\]
the real variety generated by the coefficient polynomials of $P^\lambda$.  Treating the
parameter~$s$ as an additional indeterminate we arrive at
\[
V(P^\lambda) = \SetOf{(s,t) \in \RR^{1+m}}{F_{P^\lambda,i,s}(t)=0 \; \text{for all} \; i\in[m+1] \;
  \text{and} \; s\in (0,1]} \; .
\]

\begin{lem} \label{lem:center}
  Choose color consecutive orderings of the vertices of $P^\lambda$ and $Q^\mu$. Then there is a
  lifting function $\omega$ of $P^\lambda \simpp Q^\mu = (P \times Q)^\omega$, such that the points
  in the variety $V_s({(P \times Q)^\omega})$ are exactly the points $(t,t') =
  (t_1,\dots,t_{m+n})\in \RR^{m+n}$ with $t \in V_s({P^\lambda})$ or $t' \in V_s({Q^\mu})$, that is,
  \[
  V_s({(P \times Q)^\omega}) = ( V_s({P^\lambda}) \times \RR^n ) \cup ( \RR^m \times V_s({Q^\mu}) ) \; .
  \]
\end{lem}

\begin{rem}
  The variety $V_s({P^\lambda})$ may be infinite, in general.
\end{rem}

\begin{figure}[htbp]\centering
  \begin{overpic}[width=.6\textwidth]{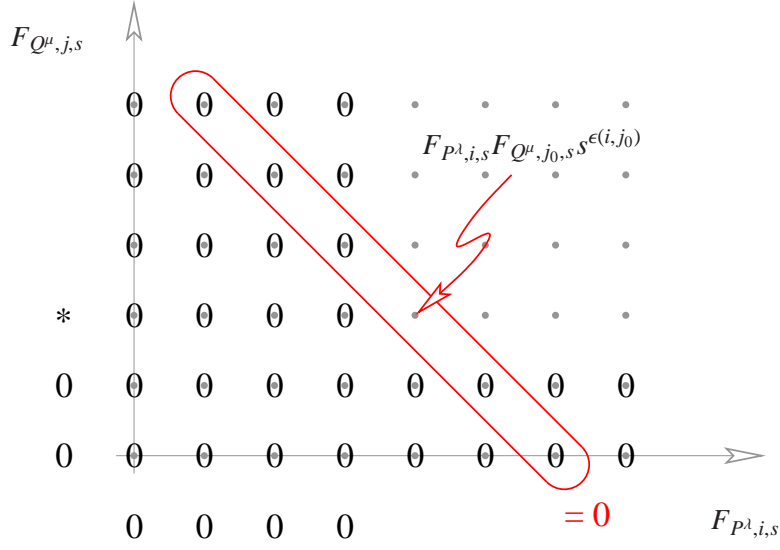}
    \put(52,54){$F_{P^\lambda,i,s} F_{Q^\mu,j_0,s} s^{\epsilon(i,j_0)}$}
    \put(92,2){$F_{P^\lambda,i,s}$}
    \put(-5,69){$F_{Q^\mu,j,s}$}
  \end{overpic}
  \caption{The inductive step in the proof of Lemma~\ref{lem:center}.  Here~$\ast$
    denotes the non-zero value of $F_{Q^\mu,j_0,s}(t')$. \label{fig:nice_product}}
\end{figure}

\begin{proof}[Proof of Lemma~\ref{lem:center}]
  For a point $t \in V_s({P^\lambda})$ we have $(t,t') \in V_s({(P \times Q)^\omega})$ for all
  $t'\in\RR^n$ by Equation~\eqref{eq:prod_coeff_poly}.  Similarly we have $(t,t') \in V_s({(P \times
    Q)^\omega})$ for $(s,t') \in V_s({Q^\mu})$ and all $t \in \RR^m$.

  For the reverse, let us assume there is a point $(t,t') \in V_s({(P \times Q)^\omega})$ but $t
  \not\in V_s({P^\lambda})$ and $t' \not\in V_s({Q^\mu})$.  Choose $i_0 \in [m+1]$ and $j_0 \in
  [n+1]$ minimal such that $F_{P^\lambda,i_0,s}(t)\not=0$ and $F_{Q^\mu,j_0,s}(t')\not=0$.  Further
  let us assume $i_0 \ge j_0$.  We prove by induction on~$i$ that $i_0>m$, or alternatively that
  $F_{P^\lambda,i,s}(t)=0$ for all $i \in [m+1]$, contradicting our assumption $t \not\in
  V_s({P^\lambda})$.
  
  We have $F_{P^\lambda,i,s}(t)=0$ for all $i<j_0$.  Note that this is also true for $j_0 =
  0$.  Now let $F_{P^\lambda,i',s}(t)=0$ for all $i'<i$.  Equation~\eqref{eq:prod_coeff_poly} yields for
  $k=i+j_0$
    \begin{align*}
      F_{(P \times Q)^\omega,i+j_0,s}(t,t') \; &= \hspace{.5cm} \sum_{i'+j'=i+j_0} F_{P^\lambda,i',s}(t) \; F_{Q^\mu,j',s}(t') \;
      s^{\epsilon(i',j')}\\
      &= \hspace{.3cm} \sum_{i'+j'=i+j_0, i'<i} F_{P^\lambda,i',s}(t) \; F_{Q^\mu,j',s}(t') \; s^{\epsilon(i',j')}\\
      & \hspace{.3cm} + \; F_{P^\lambda,i,s}(t) \; F_{Q^\mu,j_0,s}(t') \; s^{\epsilon(i,j_0)} \\
      & \hspace{.3cm} + \; \sum_{i'+j'=i+j_0, i'>i} F_{P^\lambda,i',s}(t) \; F_{Q^\mu,j',s}(t') \;
      s^{\epsilon(i',j')}\\
      &= \hspace{.2cm} 0\; ,
    \end{align*}
    since we assumed $(t,t') \in V_s({(P \times Q)^\omega})$.

    We have $F_{P^\lambda,i',s}(t) = 0$ for $i'<i$ by induction and $i'>i$ implies $j < j_0$ hence
    $F_{Q^\mu,j,s}(t') = 0$ for $i'>i$.  We are left with $F_{P^\lambda,i,s}(t) \;
    F_{Q^\mu,j_0,s}(t') \; s^{\epsilon(i,j_0)} = 0$ which in turn yields $F_{P^\lambda,i,s}(t) =0$
    since $s^{\epsilon(i,j_0)} > 0$ and $F_{Q^\mu,j_0,s}(t')\not=0$; see
    Figure~\ref{fig:nice_product}.
\end{proof}

Now we are ready to state and prove our main result.

\begin{thm}[Algebraic Product Theorem] \label{thm:alg_product_theorem}
  Let $P\subset\RR_{\ge 0}^m$ and $Q\subset\RR_{\ge 0}^n$ be non-negative full-dimensional lattice
  polytopes with {rdf}-triangulations $P^\lambda$ and $Q^\mu$ which are nice
  for some value $s_0\in(0,1]$.  Further choose any color consecutive vertex orderings for
  $P^\lambda$ and $Q^\mu$.  Then there is a lifting function $\omega: (P \times Q) \cap
  \ZZ^{m+n}\to \RR$ such that $(P \times Q)^\omega = P^\lambda \simpp Q^\mu$ is nice for~$s_0$.
  Moreover, the number of real solutions of any Wronski polynomial system associated with $(P \times
  Q)^\omega$ is bounded from below by
  \[\sigma\left((P \times Q)^\omega\right)=
  \sigma_{m,n}\;\sigma(P^\lambda)\;\sigma(Q^\mu)\;.\]
\end{thm}

\begin{proof}
  The orientability of $Y_{P \times Q}^+$ is a consequence of Corollary~\ref{cor:oriented}.  Now
  Lemma~\ref{lem:center} provides a lifting function $\omega: (P \times Q) \cap \ZZ^{m+n}\to \RR$ of
  $P^\lambda \simpp Q^\mu$ such that the $s$-degeneration $s.Y_{(P \times Q)^\omega}$ does not meet
  the center of the Wronski projection for~$s \in (0,s_0]$ and $(t,t') \in \RR^{m+n}$: Since
  $V_s({P^\lambda}) = V_s({Q^\mu}) = \emptyset$ for all $s\in (0,s_0]$ we have $V_s({(P\times
    Q)^\omega}) = ( V_s({P^\lambda}) \times \RR^n ) \cup ( \RR^m \times V_s({Q^\mu}) ) = \emptyset$
  for all $s\in (0,s_0]$.  The claim hence follows from Theorem~\ref{thm:Soprunova_Sottile} and our
  Combinatorial Product Theorem~\ref{thm:comb_product_thm}.
\end{proof}

\begin{rem} \label{rem:simpp_and_algebra}
  The decomposition $\sigma(P^\lambda\simpp Q^\mu)=\sigma_{m,n}\;\sigma(P^\lambda)\;\sigma(Q^\mu)$
  from Theorems~\ref{thm:comb_product_thm} and~\ref{thm:alg_product_theorem} reflects the geometric
  situation in the following sense: Let $M = |P\cap\ZZ^m|$ and $N = |Q\cap\ZZ^n|$ denote the number
  of lattice points of~$P$ and~$Q$, respectively.  The Wronski projection $\pi_E:
  \CC\PP^{M-1}\setminus E \to \CC\PP^m$ (and its center~$E$) depends solely on the lifting function
  $\lambda: \RR^m \to \RR$ which induces the {rdf}-triangulation~$P^\lambda$ on~$P$.
  Hence we will denote the Wronski projection~$\pi_E$ associated with~$P^\lambda$ by~$\pi_{P^\lambda}$, and its lifting to~$\Sph^{M-1}$ by~$\pi_{P^\lambda}^+$.
  To give a lower bound on the number of real roots of
  the Wronski system associated with $(P\times Q)^\omega = P^\lambda \simpp Q^\mu$ we have to bound
  the topological degree of the map $\pi_{(P\times Q)^\omega}^+$
  restricted to~$Y_{P\times Q}^+$.  A decomposition of~$\pi_{(P\times Q)^\omega}^+$ by the maps
  $\pi_{P^\lambda}^+$, $\pi_{Q^\mu}^+$,~$\pi_{\Delta_m\simpp\Delta_n}^+$, and the covers of the
  Segre embeddings is given by the following diagram which commutes provided that the lifting
  functions match as in Equation~\eqref{eq:omega2}.  Here the vertical arrows indicate the covers of
  the Segre embeddings of the appropriate dimensions.
  \[
  \xymatrix{
    \Sph^{MN-1} & 
    Y_{P\times Q}^+ \ar@{_{(}->}[l] \ar[r]^{\pi_{(P\times Q)^\omega}^+} &
    \Sph^{m+n}  & 
    Y_{\Delta_m\times\Delta_n}^+ \ar[l]_{\pi_{\stc_{m,n}}^+} \ar@{^{(}->}[r] &
    \Sph^{mn+m+n}\\
    \Sph^{M-1}\times\Sph^{N-1} \ar[u]_{\iota} &
    Y_P^+\times Y_Q^+ \ar@{_{(}->}[l] \ar[rr]^{\pi_{P^\lambda}\times\pi_{Q^\mu}}  \ar[u]_{\iota} & &
    Y_{\Delta_m}^+\times Y_{\Delta_n}^+ \ar[u]_{\iota} &
    \Sph^m\times\Sph^n \ar@{=}[l] \ar[u]_{\iota}
  }
\]
 This decomposition
 of~$\pi_{(P\times Q)^\omega}^+$ yields the
 decomposition of $\sigma(P^\lambda\simpp Q^\mu)$ given in the Theorems~\ref{thm:comb_product_thm}
 and~\ref{thm:alg_product_theorem}.
\end{rem}

\section{Cubes}
\label{sec:cubes}

We define the \emph{signature} of a lattice polytope~$P$, denoted as $\sigma(P)$, as the maximum of
the signatures of all {rdf}-triangulations of~$P$.  The signature is undefined if $P$ does not
admit any such triangulation as in Example~\ref{exmp:undef}.  However, here we are concerned with
cubes, which do have {rdf}-triangulations: This is an immediate consequence of the Product
Theorem~\ref{thm:comb_product_thm} since $C_d=[0,1]^d=I\times\dots\times I$ can be triangulated as
the $d$-fold simplicial product $I\simpp\dots\simpp I$ with zero signature.

Since $C_d$ does not contain any non-vertex lattice points, each lattice triangulation of~$C_d$ is
dense.  Note that $C_d$ does have non-regular triangulations for $d\ge 4$; see De
Loera~\cite{MR1380393}.

\subsection{Regular and Foldable Triangulations with Large Signature}
\label{sec:large_signature}

Since the simplicial product of unimodular triangulations is again unimodular it follows that each
$d$-fold simplicial product $I\simpp\dots\simpp I$ has $d!$ facets, which is the maximum that can be
obtained for the $d$-cube without introducing new vertices.  On the other hand the minimal number of
facets in a triangulation of~$C_d$ is known only for $d\le 7$; see Anderson and
Hughes~\cite{MR1411113}.  The best currently known upper and lower bounds are due to
Smith~\cite{MR1737333}, Orden and Santos~\cite{MR2013970}, and Bliss and Su~\cite{math.CO/0310142}.
For a recent survey on cubes, their triangulations, and related issues see Zong~\cite{MR2133310}.
Rambau's program \texttt{TOPCOM} allows to enumerate all regular triangulations of $C_d$ for $d\le
4$~\cite{topcom}.  This then yields the following result.

\begin{prop}\label{prop:cube:le4}
  We have $\sigma(C_1)=1$, $\sigma(C_2)=0$, $\sigma(C_3)=4$, and $\sigma(C_4)=2$.
\end{prop}

The cases of $C_1=I$ and $C_2$ are trivial.  The unique (regular and) foldable triangulation of $C_3$
with the maximal signature~$4$ is the unique minimal triangulation; it has one (black) facet of
normalized volume~$2$ and four (white) facets of normalized volume~$1$.

There is one further ingredient which relies on an explicit construction, a triangulation of $C_6$
with a non-trivial signature.  We give more details on our experiments in
Section~\ref{sec:experiments} below.

\begin{prop}\label{prop:cube:6}
  We have $\sigma(C_6) \geq 4$.
\end{prop}

\begin{thm}\label{thm:cube}
  The signature of $C_d$ for~$d \geq 3$ is bounded from below by
  \[
  \sigma(C_d) \; \geq \;
  \begin{cases}
    \; 2^{ \frac{d+1}{2} } \left( \frac{d-1}{2} \right) ! &\text{if} \; d \equiv 1 \mod 2\\
    \; \left( \frac{d}{2} \right) ! &\text{if} \; d \equiv 0 \mod 4\\
    \; \frac{2}{3} \left( \frac{d}{2} \right) !  &\text{if} \; d \equiv 2 \mod 4 \; .
  \end{cases}
  \]
\end{thm}
\begin{proof}
  Let us start with the case~$d$ odd.  Here for $C_3$ we choose the
  {rdf}-triangulation with signature~$4$ from Proposition~\ref{prop:cube:le4}.  For $d\ge 5$ we factorize
  $C_d$ as $C_2\times C_{d-2}$ and choose a color consecutive vertex ordering for~$C_{d-2}$.  There
  is only one triangulation to choose for~$C_2$, but we take the symmetric ordering of the vertices;
  see Example~\ref{exmp:bipyramid}.  The signature of $stc_{2,d-2}$ equals $(d-1)/2$ by
  Proposition~\ref{prop:sign_m_n} and the second case of Proposition~\ref{prop:sign_of_bipyramid}
  inductively gives
  \[
  \sigma(C_d) \;\ge\; 2\:\sigma_{d-2,2}\:\sigma(C_{d-2})\;\ge\;
  2\:\frac{d-1}{2}\:2^{\frac{d-3}{2}}\:\left(\frac{d-3}{2}\right)!
  \;=\;2^{ \frac{d+1}{2} }\:\left( \frac{d-1}{2} \right) ! \;.
  \]
  
  If $d\equiv 0\mod 4$ then we inductively prove that $\sigma(C_d)\ge \left( \frac{d}{2} \right) !$.
  The induction starts with $d=4$ by Proposition~\ref{prop:cube:le4}. For~$d \ge 8$ we decompose
  $C_d$ as $C_4\times C_{d-4}$. The signature of~$\stc_{4,d-4}$ equals $d(d-2)/8$ by
  Proposition~\ref{prop:sign_m_n}. Choosing color consecutive orderings for $C_4$ and $C_{d-4}$
  Theorem~\ref{thm:comb_product_thm} now yields
  \[
  \sigma(C_d) \;\ge\; \sigma_{4,d-4}\:\sigma(C_4)\:\sigma(C_{d-4})
  \;\ge\; \frac{d(d-2)}{8}\:2\:\left(\frac{d-4}{2}\right)! 
  \;=\; \left(\frac{d}{2}\right)! \;.
  \]
  
  In the remaining case where $d \equiv 2 \mod 4$ we construct~$C_d$ as a simplicial product
  of~$C_6$ and~$C_{d-6}$.  By the explicit construction in Proposition~\ref{prop:cube:6} the
  signature of~$C_6$ is at least~$4$.  The signature of~$C_{d-6}$ is bounded from below by
  $(d-6)/2!$ as just proved.  Proposition~\ref{prop:sign_m_n} yields $\sigma_{6,d-6}={\frac{d}{2}
    \choose 3}$, and Theorem~\ref{thm:comb_product_thm} completes the proof:
  \[
  \sigma(C_d) \;\ge\; \sigma_{6,d-6}\:\sigma(C_6)\:\sigma(C_{d-6}) \;\ge\;
  \frac{ \frac{d}{2} \left( \frac{d}{2}-1 \right) \left( \frac{d}{2}-2 \right) }{3!}\:4\:
  \left( \frac{d}{2}-3 \right) !
  \;=\; \frac{2}{3}\:\left( \frac{d}{2} \right) ! \;.
  \]
\end{proof}

\subsection{Nice Triangulations}
Our main result, the Algebraic Product Theorem~\ref{thm:alg_product_theorem}, asserts that the
simplicial product of two nice triangulations~$P^\lambda$ and~$Q^\mu$ is again nice, provided that
the vertex ordering of~$P^\lambda$ and~$Q^\mu$ are color consecutive.  So what about the
triangulations of the $d$-cube with signature in $\Omega(\lceil d/2 \rceil !)$ constructed in
Section~\ref{sec:large_signature} above?  Since the construction for~$d$ odd was based on the
symmetric vertex ordering for the square, which is not color consecutive,
Theorem~\ref{thm:alg_product_theorem} does not apply.  The goal of this section is thus to construct
nice cube {rdf}-triangulations from a decomposition into different factors.

The \emph{geometric signature} $\sigma^+(P)$ of a lattice polytope~$P$ is defined as the maximum of
the signatures of all {rdf}-triangulations of~$P$ which are nice for some parameter $s\in(0,1]$.
Clearly, $\sigma^+(P)\le\sigma(P)$.  Note that $Y_{C_d}^+$ is always oriented by
Corollary~\ref{cor:oriented} since $C_d=I\times I\times \dots \times I$, and $I$ is Cox-oriented.

Let us examine two cases of low dimension explicitly: There is a lifting function $C_3 \cap \ZZ^3
\to \NN$ such that the induced triangulation is the unique minimal triangulation of the $3$-cube
from Proposition~\ref{prop:cube:le4}, and the toric degeneration meets the center only for $s=1$;
see~\cite{math.AG/0409504}.  This implies $\sigma^+(C_3)=4$.  In the subsequent
Section~\ref{sec:experiments} a triangulation~$C_4^\lambda$ of the $4$-cube with signature equal
to~$2$ is constructed explicitly via a lifting function $\lambda: C_4 \cap \ZZ^4 \to \NN$.  The
variety~$V({C_4^\lambda})$ (see Section~\ref{sec:poly_sys_from_prod}), describing the values of~$s$
for which the center of the projection is met, consists of two isolated points for some $s_1>1$ and
some $s_2<0$, hence~$C_4^\lambda$ is nice for any $s_0 \in (0,1]$.  A complete enumeration of all
regular triangulation of $C_4$ shows that $\sigma^+(C_4)=2$.

We want to avoid to split off factors which are squares, since neither of its two vertex orderings
can be used for our purposes: The color consecutive vertex ordering has signature zero, and products
with respect to the symmetric vertex ordering are not known to be nice.  Hence we factorize
\[
C_d \;=\; \begin{cases}
  C_1 \times C_{d-1} &\text{if} \; d \equiv 1 \mod 4\\
  C_3 \times C_{d-3} &\text{if} \; d \equiv 3 \mod 4 \; ,
\end{cases}
\]
which means that we reduced the cases $d\equiv 1\mod 4$ and $d\equiv 3\mod 4$ to the case
\mbox{$d\equiv 0\mod 4$}.  Proposition~\ref{prop:sign_m_n} and Theorem~\ref{thm:comb_product_thm}
yield for $d\equiv 1\mod 4$
\[
\sigma^+(C_d) \;\ge\; \sigma_{1,d-1}\:\sigma^+(C_1)\:\sigma^+(C_{d-1}) \;=\; \sigma^+(C_{d-1})
\;\ge\; \left( \frac{d-1}{2} \right) ! \;.
\]
For $d\equiv 3\mod 4$ we have
\[
\sigma^+(C_d) \;\ge\; \sigma_{3,d-3}\:\sigma^+(C_3)\:\sigma^+(C_{d-3}) \;\ge\; 
\frac{d-1}{2}\:4\:\left( \frac{d-3}{2} \right) ! \;=\;
4\:\left( \frac{d-1}{2} \right) ! \;,
\]
and we obtain an overall lower bound in $\Omega(\lfloor d/2 \rfloor!)$ for the geometric signature
of the $d$-cube.  Observe that this lower bound for the signature in the case of~$d$ odd is significantly weaker than
the bound given in Theorem~\ref{thm:cube}, which does not take the geometric properties of the
Wronski projection into account.


\begin{cor}\label{cor:nice_cubes}
  For $3 \leq d\not\equiv 2\mod 4$ there are rdf-triangulations of the $d$-cube with signature at least~$\lfloor d/2 \rfloor !$ which are nice for any $s_0 \in (0,1)$.
\end{cor}

Proving that the triangulation of the $6$-cube with signature~$4$ from Proposition~\ref{prop:cube:6}
(together with its generating lifting function) is nice for some $s_0 \in (0,1]$ would also settle
the $d\equiv 2\mod 4$ case.  However, with the techniques of Section~\ref{sec:experiments} one needs
to solve a system of seven polynomials in the seven unknowns $s, x_1, \dots, x_6$ of maximal total
degree~$386$; see Problem~\ref{prob:C6}.  This is beyond the scope of this paper.

\subsection{Constructions and Computer Experiments}\label{sec:experiments}

We completely enumerated all regular triangulations of the $d$-cube~$C_4$ up to symmetry using
\texttt{TOPCOM}~\cite{topcom}.  These $235,\!277$
triangulations were then checked whether they are foldable by
\texttt{polymake}~\cite{polymake,DMV:polymake,math.CO/0507273}; it turns out that their total number
is~$454$.  For all the foldable ones we computed the signature, and we found~$36$ triangulations
with signature~$2$, all other foldable triangulations of $C_4$ have a vanishing signature.  The
regularity of Example~\ref{exmp:4cube:signature2} was independently verified by the explicit
construction of a lifting function.

\begin{exmp}\label{exmp:4cube:signature2}
  We now give an explicit description of an {rdf}-triangulation $C_4^\lambda$ of the
  $4$-cube with signature two.  To this end we encode the vertices of $C_4$, that is, the
  $0/1$-vectors of length~$4$ as the hexadecimal digits $0,1,\dots,9,a,b,c,d,e,f$.  The lifting
  function $\lambda$ and the vertex $5$-coloring is given in Table~\ref{tab:4cube:col_lift}.
  The facets of $C_4^\lambda$ are listed in Table~\ref{tab:4cube:facets}, and the $f$-vector reads
  $(16,64,107,81,23)$.
  
  \begin{table}[htbp]\centering
    \caption{The vertex $5$-coloring $c$ and a lifting function $\lambda$ for $C_4^\lambda$
      described in Example~\ref{exmp:4cube:signature2}.  The vertices of the
      first facet ${0 1 2 4 8}$ are chosen as the colors.\label{tab:4cube:col_lift}}
    \renewcommand{\arraystretch}{0.9}
    \begin{tabular*}{.9\linewidth}{@{\extracolsep{\fill}}lrrrrrrrrrrrrrrrr@{}}\toprule      
      $v$ & $0$ & $1$ & $2$ & $3$ & $4$ & $5$ & $6$ & $7$ & $8$ & $9$ & $a$ & $b$ & $c$ & $d$
      & $e$ & $f$\\
      \midrule
      $\lambda(v)$ & $0$ & $0$ & $0$ & $4$ & $0$ & $2$ & $8$ & $8$ & $10$ & $11$ & $19$ & $19$ & $10$ & $19$ & $24$ & $31$\\
      $c(v)$ & $0$ & $1$ & $2$ & $4$ & $4$ & $0$ & $0$ & $1$ & $8$ & $2$ & $1$ & $0$ & $2$ & $4$ & $4$ & $8$\\
      \bottomrule
    \end{tabular*}
  \end{table}

  \begin{table}[htbp]\centering
    \caption{Facets of the triangulation $C_4^\lambda$.\label{tab:4cube:facets}}
    \renewcommand{\arraystretch}{0.9}
    \begin{tabular*}{.9\linewidth}{@{\extracolsep{\fill}}cccccccc@{}}\toprule      
      ${0 1 2 4 8}$ &
      ${1 2 3 5 8}$ &
      ${1 2 4 5 8}$ &
      ${1 3 5 8 9}$ &
      ${2 3 7 8 b}$ &
      ${2 3 5 7 8}$ &
      ${2 4 5 7 8}$ &
      ${2 4 6 7 8}$ \\
      ${2 6 7 8 e}$ &
      ${2 7 8 b e}$ &
      ${2 8 a b e}$ &
      ${3 5 7 8 9}$ &
      ${3 7 8 9 b}$ &
      ${4 5 7 8 c}$ &
      ${4 6 7 8 c}$ &
      ${5 7 8 9 d}$ \\
      ${5 7 8 c d}$ &
      ${6 7 8 c e}$ &
      ${7 8 9 b d}$ &
      ${7 8 b c d}$ &
      ${7 8 b c e}$ &
      ${7 b c e f}$ &
      ${7 b c d f}$ \\
      \bottomrule
    \end{tabular*}
  \end{table}

  As mentioned before, the double cover~$Y_{C_d}^+$ of the associated real
  toric variety of the $d$-cube is indeed oriented for all dimensions~$d$.
  To prove that~$C_4^\lambda$ is nice for any $s_0 \in (0,1]$ we examine the
  variety~$V({C_4^\lambda})$, describing the values of~$s$ for which the center of the projection is
  met; see Section~\ref{sec:poly_sys_from_prod}.  The variety~$V({C_4^\lambda})$ is the solution set
  of the ideal~$I({C_4^\lambda})$ generated by the five coefficient polynomials
   \begin{align*}
     &F_{C_4^\lambda,0,s} = 1 + s^2 x_1 x_3 + s^8 x_2 x_3 + s^{19} x_1 x_2 x_4\; ,\\
     &F_{C_4^\lambda,1,s} = x_1 + s^8 x_1 x_2 x_3 + s^{19} x_2 x_4\; ,\\
     &F_{C_4^\lambda,2,s} = x_2 + s^{10} x_3 x_4 + s^{11} x_1 x_4\; ,\\
     &F_{C_4^\lambda,3,s} = x_3 + s^4 x_1 x_2 + s^{19} x_1 x_3 x_4 + s^{24} x_2 x_3 x_4\; , \; \text{and}\\
     &F_{C_4^\lambda,4,s} = x_4 + s^{31} x_1 x_2 x_3 x_4\; .
   \end{align*}

   For the lexicographical ordering $x_4>x_3>x_2>x_1>s$ a Gr\"obner basis of~$I({C_4^\lambda})$ reads
   (computed by
   \texttt{MAGMA}~\cite{magma})
   \[
   \left\{\: x_4 + g_4(s),\: x_3 + g_3(s),\: x_2 + g_2(s),\: x_1 + g_1(s),\: g_s(s) \:\right\}\; ,
   \]
   for certain polynomials $g_s,g_1,\dots,g_4 \in \QQ[s]$.  The polynomial $g_s(s)$ is displayed in
   Figure~\ref{fig:gs}, and the others are by far too large to be listed.  The essential feature of
   this Gr\"obner basis is that knowing the (real) roots of the polynomial $g_s(s)$ of degree~$444$
   allows to compute the values for $x_1,\dots,x_4$ directly.

\begin{table}[bhtp]\centering
  \caption{Approximate coordinates for the two points in the variety $V({C_4^\lambda})$.\label{tab:4cube:approx}}
  \renewcommand{\arraystretch}{0.9}
  \begin{tabular*}{.9\linewidth}{@{\extracolsep{\fill}}lr@{.\extracolsep{0pt}}l@{\extracolsep{\fill}}r@{.\extracolsep{0pt}}l@{}}\toprule
    $s$   & -0&9955941875452 &  1&0003839818262 \\
    $x_1$ &  1&3469081499925 & -1&1340421741317 \\
    $x_2$ &  0&7663015145691 & -1&8447577233888 \\
    $x_3$ &  1&1109881050869 & -0&4723488390037 \\
    $x_4$ &  3&4823714929884 & -1&1436761629897 \\
    \bottomrule
  \end{tabular*}      
\end{table}  

   It turns out that $g_s(s)$ has exactly two real roots $s_1$ and $s_2$ with $s_1>1$ and
   $-1<s_2<0$.  Given $g_s(s)$ this can be verified with any standard computer algebra program by
   computing all~$444$ distinct (complex) solutions.  Additionally, this was counter-checked via
   Collins' method of cylindrical algebraic decomposition~\cite{MR0403962}, as implemented in
   \texttt{QEPCAD}~\cite{qepcad}.  Approximate values for the two real zeroes of $g_s$ are given in
   Table~\ref{tab:4cube:approx}.  It follows that~$C_4^\lambda$ is nice for any $s_0 \in (0,1]$.
 \end{exmp}
  
\begin{figure}[t]
  \begin{minipage}{.87\textwidth}
    \include{gs}
  \end{minipage}
  \caption{The polynomial $g_s(s)$ of the Gr\"obner basis of~$I({C_4^\lambda})$.}
  \label{fig:gs}
\end{figure}
 
\enlargethispage{\baselineskip}
While, with current computers, it seems to be out of reach to completely enumerate all
triangulations of most polytopes in dimension~$5$ and beyond, \texttt{TOPCOM} can still be used to
enumerate large numbers of triangulations.  We let \texttt{TOPCOM} compute altogether $59,\!083$
different triangulations which originate from randomly chosen placing triangulations by successive
flipping.  Not a single triangulation among these was foldable.  Next we took the triangulation of
$C_5$ with signature~$16$ that comes from Theorem~\ref{thm:cube} and we inspected $102,\!184$
triangulations by random flipping.  This way we found only two more foldable triangulations, one
with signature~$14$ and a second one with signature~$16$.
 
For $C_6$ the situation is more complicated.  None of our results so far directly yields any
foldable triangulation with a positive signature: All the simplicial product triangulations of $C_6$
arising from decomposing $C_6$ as a product of two (or more) cubes of smaller dimensions do not
yield a non-trivial lower bound since at least one factor vanishes in the corresponding expressions
in Proposition~\ref{prop:sign_of_bipyramid} and Theorem~\ref{thm:comb_product_thm}.  And, as can be
expected from the $5$-dimensional case, \texttt{TOPCOM} did not find a foldable triangulation with a
positive signature either.  Therefore we took a detour in that we used \texttt{TOPCOM} to study
triangulations of the product of the $4$-simplex and the square.  This time we were lucky to find a
foldable triangulation with signature~$2$, which also turned out to be regular.

\begin{prop}\label{prop:delta4xc2}
  We have $\sigma(\Delta_4\times C_2)\ge 2$
\end{prop}
 
In the sequel we denote this {rdf}-triangulation of $\Delta_4 \times C_2$ with signature~$2$ by $S$,
and let $C_4^\lambda$ be the {rdf}-triangulation of $C_4$ with signature~$2$ from
Proposition~\ref{prop:cube:le4}.  Then the product $C_6=C_4\times C_2$ inherits a polytopal
subdivision into facets of type $\Delta_4 \times C_2$ from $C_4^\lambda$.  Each of these facets can
now be triangulated using~$S$ in such a way that one obtains an {rdf}-triangulation of $C_6$ with
signature~$4$.  Its $f$-vector equals $(64,656,2640,5298,5676,3115,690)$.  This establishes
Proposition~\ref{prop:cube:6}.

\begin{prob}\label{prob:C6}
  In order to decide whether the triangulation of~$C_6$ from Proposition~\ref{prop:cube:6} (together
  with its generating lifting function) is nice for some $s_0 \in (0,1]$, it suffices to prove that
  the real variety generated by
   \begin{align*}
     F_{C_6,0,s} = \;&1+s^2x_5x_6+s^8x_1x_6+s^{55}x_1x_3+s^{57}x_1x_3x_5x_6+s^{124}x_2x_3+s^{151}x_2x_3x_5x_6+\\
                   &s^{157}x_1x_2x_3x_6+s^{197}x_1x_2x_4+s^{218}x_2x_4x_6+s^{224}x_1x_2x_4x_5x_6,\\
     F_{C_6,1,s} = \;&x_6+s^4x_1x_5+s^{41}x_2x_5x_6+s^{55}x_1x_3x_6+s^{122}x_1x_4x_5x_6+s^{128}x_1x_2x_3x_5+\\
                   &s^{149}x_2x_3x_6+s^{167}x_3x_4x_5x_6+s^{189}x_2x_4x_5+s^{222}x_1x_2x_4x_6,\\
     F_{C_6,2,s} = \;&x_5+s^8x_1x_5x_6+s^{55}x_1x_3x_5+s^{124}x_2x_3x_5+s^{157}x_1x_2x_3x_5x_6+s^{197}x_1x_2x_4x_5+\\
                   &s^{218}x_2x_4x_5x_6,\\
     F_{C_6,3,s} = \;&x_1+s^8x_2x_5+s^{35}x_3x_6+s^{55}x_4x_5x_6+s^{89}x_1x_4x_5+s^{92}x_1x_2x_6+s^{124}x_1x_2x_3+\\
                   &s^{134}x_3x_4x_5+s^{185}x_2x_4+s^{218}x_1x_3x_4x_6+s^{311}x_2x_3x_4x_6+s^{380}x_1x_2x_3x_4x_5x_6,\\
     F_{C_6,4,s} = \;&x_2+s^{10}x_3x_5+s^{39}x_4x_6+s^{67}x_1x_2x_5+s^{81}x_1x_4+s^{126}x_3x_4+s^{193}x_1x_3x_4x_5+\\
                   &s^{286}x_2x_3x_4x_5+s^{364}x_1x_2x_3x_4x_6,\\
     F_{C_6,5,s} = \;&x_3+s^{12}x_4x_5+s^{37}x_2x_6+s^{57}x_1x_2+s^{118}x_1x_4x_6+s^{163}x_3x_4x_6+s^{183}x_1x_3x_4+\\
                   &s^{276}x_2x_3x_4+s^{337}x_1x_2x_3x_4x_5, \; \text{and}\\
     F_{C_6,6,s} = \;&x_4+s^{49}x_3x_5x_6+s^{106}x_1x_2x_5x_6+s^{325}x_1x_2x_3x_4+s^{325}x_2x_3x_4x_5x_6+\\
                   &s^{232}x_1x_3x_4x_5x_6
   \end{align*}
   is empty for all $s\in (0,s_0]$.  We leave this as an open problem.
\end{prob}

\section{A Further Remark and Several Acknowledgments}

\enlargethispage{\baselineskip}
Triangulations of the rectangular grid~$G_{k,l} = [0,k] \times [0,l]$ are an interesting subject of
its own; see, for instance, Kaibel and Ziegler~\cite{MR2011739} and the references there.  Note that
each triangulation of the grid is dense if and only if it is unimodular.  Even without the
assumption of regularity we do not know of a single dense and foldable triangulation of $G_{k,l}$
with a positive signature.

\begin{prob}
  For which parameters $k$ and $l$, if any, does the rectangular grid $G_{k,l}$ admit a unimodular
  and foldable triangulation with a positive signature?
\end{prob}

Till Stegers helped with Gr\"obner bases computations.  Chris W. Brown gave a first computer based
proof of the fact that the variety $V({C_4^\lambda})$ consists of two isolated points via
\texttt{QEPCAD}~\cite{qepcad}, and he also provided our approximate coordinates. Frank Sottile
helped us to a better understanding of the geometric situation. In particular, he noticed that the
product of two Cox-orientable lattice polytopes is not necessarily Cox-orientable.  Two referees
gave very useful comments on a previous version.  We are indebted to all of them.  Finally, we are
grateful to Thorsten Theobald for stimulating discussions on the subject.

\bibliographystyle{amsplain}
\bibliography{main}

\end{document}